\magnification=\magstep1
\hsize16truecm
\vsize23.5truecm
\topskip=1truecm
\raggedbottom
\abovedisplayskip=3mm
\belowdisplayskip=3mm
\abovedisplayshortskip=0mm
\belowdisplayshortskip=2mm
\normalbaselineskip=12pt
\normalbaselines
\font\titlefont= cmcsc10 at 12pt
\def\chno{1}
\def\F{\Bbb F}

\def\Q{\Bbb Q}
\def\P{\Bbb P}

\def\G{\Bbb G}

\def\rk{\rm rank}
 
%
\catcode`\@=11
\font\tenmsa=msam10
\font\sevenmsa=msam7
\font\fivemsa=msam5
\font\tenmsb=msbm10
\font\sevenmsb=msbm7
\font\fivemsb=msbm5
\newfam\msafam
\newfam\msbfam
\textfont\msafam=\tenmsa  \scriptfont\msafam=\sevenmsa
  \scriptscriptfont\msafam=\fivemsa
\textfont\msbfam=\tenmsb  \scriptfont\msbfam=\sevenmsb
  \scriptscriptfont\msbfam=\fivemsb
\def\hexnumber@#1{\ifcase#1 0\or1\or2\or3\or4\or5\or6\or7\or8\or9\or
	A\or B\or C\or D\or E\or F\fi }
\edef\msa@{\hexnumber@\msafam}
\edef\msb@{\hexnumber@\msbfam}
\mathchardef\square="0\msa@03
\mathchardef\subsetneq="3\msb@28
\mathchardef\ltimes="2\msb@6E
\mathchardef\rtimes="2\msb@6F
\def\Bbb{\ifmmode\let\next\Bbb@\else
	\def\next{\errmessage{Use \string\Bbb\space only in math mode}}\fi\next}
\def\Bbb@#1{{\Bbb@@{#1}}}
\def\Bbb@@#1{\fam\msbfam#1}
\catcode`\@=12
%
%
%
\def\mapdown#1{\Big\downarrow\rlap{$\vcenter{\hbox{$\scriptstyle#1$}}$}}

\def\mapright#1{\ \smash{\mathop{\longrightarrow}\limits^{#1}}\ }

\def\tto{\longrightarrow}

\def\lllongrightarrow{\relbar\joinrel\relbar\joinrel\relbar\joinrel\rightarrow}
\def\ttto{\lllongrightarrow}
\input amssym.def
\input amssym.tex
\def\today{\ifcase\month\or
 Jan\or Febr\or  Mar\or  Apr\or May\or Jun\or  Jul\or
 Aug\or  Sep\or  Oct\or Nov\or  Dec\or\fi
 \space\number\day, \number\year}
\noindent
\font\eighteenbf=cmbx10 scaled\magstep3
\font\titlefont=cmcsc10 at 12pt
\vskip 5.0pc
\centerline{\eighteenbf On a  Stratification of the Moduli} 
\bigskip
\centerline{\eighteenbf of K3 Surfaces }
\vskip 2pc
\centerline{\titlefont G.\ van der Geer and T.\ Katsura}
\vskip2pc
\catcode`@=11
\message{Loading the plain-augmented-CM format}%
%
%
%
\font\eighteenrm=cmr10 scaled\magstep3          
\font\fourteenrm=cmr10 scaled\magstep2          
\font\twelverm=cmr12                            
\font\elevenrm=cmr10 scaled\magstephalf         
\font\tenrm=cmr10
\font\ninerm=cmr9
\font\eightrm=cmr8
\font\sevenrm=cmr7
\font\sixrm=cmr6
\font\fiverm=cmr5
%
\font\eighteeni=cmmi10 scaled\magstep3          
\font\fourteeni=cmmi10 scaled\magstep2          
\font\twelvei=cmmi12                            
\font\eleveni=cmmi10 scaled\magstephalf         
\font\teni=cmmi10
\font\ninei=cmmi9
\font\eighti=cmmi8
\font\seveni=cmmi7
\font\sixi=cmmi6
\font\fivei=cmmi5
%
\font\eighteensy=cmsy10 scaled\magstep3         
\font\fourteensy=cmsy10 scaled\magstep2         
\font\twelvesy=cmsy10 scaled\magstep1           
\font\elevensy=cmsy10 scaled\magstephalf        
\font\tensy=cmsy10
\font\ninesy=cmsy9
\font\eightsy=cmsy8
\font\sevensy=cmsy7
\font\sixsy=cmsy6
\font\fivesy=cmsy5
%
\font\eighteenex=cmex10 scaled\magstep3         
\font\fourteenex=cmex10 scaled\magstep2         
\font\twelveex=cmex10 scaled\magstep1           
\font\elevenex=cmex10 scaled\magstephalf        
\font\tenex=cmex10
%
\font\eighteenbf=cmbx10 scaled\magstep3         
\font\fourteenbf=cmbx10 scaled\magstep2         
\font\twelvebf=cmbx12                           
\font\elevenbf=cmbx10 scaled\magstephalf        
\font\tenbf=cmbx10
\font\ninebf=cmbx9
\font\eightbf=cmbx8
\font\sevenbf=cmbx7
\font\sixbf=cmbx6
\font\fivebf=cmbx5
%
\font\eighteentt=cmtt10 scaled\magstep3         
\font\fourteentt=cmtt10 scaled\magstep2         
\font\twelvett=cmtt12                           
\font\eleventt=cmtt10 scaled \magstephalf       
\font\tentt=cmtt10
\font\ninett=cmtt9
\font\eighttt=cmtt8
%
\font\eighteensl=cmsl10 scaled\magstep3         
\font\fourteensl=cmsl10 scaled\magstep2         
\font\twelvesl=cmsl12                           
\font\elevensl=cmsl10 scaled \magstephalf       
\font\tensl=cmsl10
\font\ninesl=cmsl9
\font\eightsl=cmsl8
%
\font\eighteenit=cmti10 scaled\magstep3         
\font\fourteenit=cmti10 scaled\magstep2         
\font\twelveit=cmti12                           
\font\elevenit=cmti10 scaled \magstephalf       
\font\tenit=cmti10
\font\nineit=cmti9
\font\eightit=cmti8
%
\font\eighteencsc=cmcsc10 scaled\magstep3         
\font\fourteencsc=cmcsc10 scaled\magstep2         
\font\twelvecsc=cmcsc10 scaled \magstep1          
\font\elevencsc=cmcsc10 scaled \magstephalf       
\font\tencsc=cmcsc10
\font\ninecsc=cmcsc9
\font\eightcsc=cmcsc8
%
%
%
%
%
%
%
\skewchar\eighteeni='177 \skewchar\fourteeni='177 \skewchar\twelvei='177
\skewchar\eleveni='177 \skewchar\teni='177 \skewchar\ninei='177
\skewchar\eighti='177 \skewchar\seveni='177 \skewchar\sixi='177
\skewchar\fivei='177
\skewchar\eighteensy='60 \skewchar\fourteensy='60 \skewchar\twelvesy='60
\skewchar\elevensy='60 \skewchar\tensy='60 \skewchar\ninesy='60
\skewchar\eightsy='60 \skewchar\sevensy='60 \skewchar\sixsy='60
\skewchar\fivesy='60
%
%
%
%
\let\usereighteenpointmacro=\relax
\def\eighteenpoint{\let\pointsize=\eighteenpoint
  \textfont0=\eighteenrm \scriptfont0=\twelverm \scriptscriptfont0=\ninerm
  \def\rm{\fam\z@\eighteenrm}%
  \textfont1=\eighteeni \scriptfont1=\twelvei \scriptscriptfont1=\ninei
  \def\mit{\fam\@ne}\def\oldstyle{\fam\@ne\eighteeni}%
  \textfont2=\eighteensy \scriptfont2=\twelvesy \scriptscriptfont2=\ninesy
  \def\cal{\fam\tw@}%
  \textfont3=\eighteenex \scriptfont3=\eighteenex \scriptscriptfont3=\eighteenex
  \textfont\itfam=\eighteenit
  \def\it{\fam\itfam\eighteenit}%
  \textfont\slfam=\eighteensl
  \def\sl{\fam\slfam\eighteensl}%
  \textfont\bffam=\eighteenbf \scriptfont\bffam=\twelvebf
  \scriptscriptfont\bffam=\ninebf
  \def\bf{\fam\bffam\eighteenbf}%
  \textfont\ttfam=\eighteentt
  \def\tt{\fam\ttfam\eighteentt}%
  \textfont\cscfam=\eighteencsc
  \def\csc{\fam\cscfam=\eighteencsc}%
  \def\big##1{{\hbox{$\left##1\vbox to15.3\p@{}\right.\n@space$}}}
  \def\Big##1{{\hbox{$\left##1\vbox to20.7\p@{}\right.\n@space$}}}
  \def\bigg##1{{\hbox{$\left##1\vbox to26.1\p@{}\right.\n@space$}}}
  \def\Bigg##1{{\hbox{$\left##1\vbox to31.5\p@{}\right.\n@space$}}}
  \setnormalsp@cing
  \usereighteenpointmacro
  }%
\let\usersixteenpointmacro=\relax
\def\sixteenpoint{%
  \message{You have asked for 16pt, but you are getting 14pt fonts.
           There are no 16pt CM fonts in the local font library.}%
  \fourteenpoint \let\pointsize=\sixteenpoint
  \usersixteenpointmacro
  }%
\let\userfourteenpointmacro=\relax
\def\fourteenpoint{\let\pointsize=\fourteenpoint
  \textfont0=\fourteenrm \scriptfont0=\tenrm \scriptscriptfont0=\sevenrm
  \def\rm{\fam\z@\fourteenrm}%
  \textfont1=\fourteeni \scriptfont1=\teni \scriptscriptfont1=\seveni
  \def\mit{\fam\@ne}\def\oldstyle{\fam\@ne\fourteeni}%
  \textfont2=\fourteensy \scriptfont2=\tensy \scriptscriptfont2=\sevensy
  \def\cal{\fam\tw@}%
  \textfont3=\fourteenex \scriptfont3=\fourteenex \scriptscriptfont3=\fourteenex
  \textfont\itfam=\fourteenit
  \def\it{\fam\itfam\fourteenit}%
  \textfont\slfam=\fourteensl
  \def\sl{\fam\slfam\fourteensl}%
  \textfont\bffam=\fourteenbf \scriptfont\bffam=\tenbf
  \scriptscriptfont\bffam=\sixbf
  \def\bf{\fam\bffam\fourteenbf}%
  \textfont\ttfam=\fourteentt
  \def\tt{\fam\ttfam\fourteentt}%
  \textfont\cscfam=\fourteencsc
  \def\csc{\fam\cscfam\fourteencsc}%
  \def\big##1{{\hbox{$\left##1\vbox to11.9\p@{}\right.\n@space$}}}
  \def\Big##1{{\hbox{$\left##1\vbox to16.1\p@{}\right.\n@space$}}}
  \def\bigg##1{{\hbox{$\left##1\vbox to20.3\p@{}\right.\n@space$}}}
  \def\Bigg##1{{\hbox{$\left##1\vbox to24.5\p@{}\right.\n@space$}}}
  \setnormalsp@cing
  \userfourteenpointmacro
  }%
\let\usertwelvepointmacro=\relax
\def\twelvepoint{\let\pointsize=\twelvepoint
  \textfont0=\twelverm \scriptfont0=\ninerm \scriptscriptfont0=\sixrm
  \def\rm{\fam\z@\twelverm}%
  \textfont1=\twelvei \scriptfont1=\ninei \scriptscriptfont1=\sixi
  \def\mit{\fam\@ne}\def\oldstyle{\fam\@ne\twelvei}%
  \textfont2=\twelvesy \scriptfont2=\ninesy \scriptscriptfont2=\sixsy
  \def\cal{\fam\tw@}%
  \textfont3=\twelveex \scriptfont3=\twelveex \scriptscriptfont3=\twelveex
  \textfont\itfam=\twelveit
  \def\it{\fam\itfam\twelveit}%
  \textfont\slfam=\twelvesl
  \def\sl{\fam\slfam\twelvesl}%
  \textfont\bffam=\twelvebf \scriptfont\bffam=\ninebf
  \scriptscriptfont\bffam=\sixbf
  \def\bf{\fam\bffam\twelvebf}%
  \textfont\ttfam=\twelvett
  \def\tt{\fam\ttfam\twelvett}%
  \textfont\cscfam=\twelvecsc
  \def\csc{\fam\cscfam\twelvecsc}%
  \def\big##1{{\hbox{$\left##1\vbox to10.2\p@{}\right.\n@space$}}}
  \def\Big##1{{\hbox{$\left##1\vbox to13.8\p@{}\right.\n@space$}}}
  \def\bigg##1{{\hbox{$\left##1\vbox to17.4\p@{}\right.\n@space$}}}
  \def\Bigg##1{{\hbox{$\left##1\vbox to21\p@{}\right.\n@space$}}}
  \setnormalsp@cing
  \usertwelvepointmacro
  }%
\let\userelevenpointmacro=\relax
\def\elevenpoint{\let\pointsize=\elevenpoint
  \textfont0=\elevenrm \scriptfont0=\eightrm \scriptscriptfont0=\fiverm
  \def\rm{\fam\z@\elevenrm}%
  \textfont1=\eleveni \scriptfont1=\eighti \scriptscriptfont1=\fivei
  \def\mit{\fam\@ne}\def\oldstyle{\fam\@ne\eleveni}%
  \textfont2=\elevensy \scriptfont2=\eightsy \scriptscriptfont2=\fivesy
  \def\cal{\fam\tw@}%
  \textfont3=\elevenex \scriptfont3=\elevenex \scriptscriptfont3=\elevenex
  \textfont\itfam=\elevenit
  \def\it{\fam\itfam\elevenit}%
  \textfont\slfam=\elevensl
  \def\sl{\fam\slfam\elevensl}%
  \textfont\bffam=\elevenbf \scriptfont\bffam=\eightbf
  \scriptscriptfont\bffam=\fivebf
  \def\bf{\fam\bffam\elevenbf}%
  \textfont\ttfam=\eleventt
  \def\tt{\fam\ttfam\eleventt}%
  \textfont\cscfam=\elevencsc
  \def\csc{\fam\cscfam\elevencsc}%
  \def\big##1{{\hbox{$\left##1\vbox to9.3\p@{}\right.\n@space$}}}
  \def\Big##1{{\hbox{$\left##1\vbox to12.6\p@{}\right.\n@space$}}}
  \def\bigg##1{{\hbox{$\left##1\vbox to16\p@{}\right.\n@space$}}}
  \def\Bigg##1{{\hbox{$\left##1\vbox to19.2\p@{}\right.\n@space$}}}
  \setnormalsp@cing
  \userelevenpointmacro
  }%
\let\usertenpointmacro=\relax
\def\tenpoint{\let\pointsize=\tenpoint
  \textfont0=\tenrm \scriptfont0=\sevenrm \scriptscriptfont0=\fiverm
  \def\rm{\fam\z@\tenrm}%
  \textfont1=\teni \scriptfont1=\seveni \scriptscriptfont1=\fivei
  \def\mit{\fam\@ne}\def\oldstyle{\fam\@ne\teni}%
  \textfont2=\tensy \scriptfont2=\sevensy \scriptscriptfont2=\fivesy
  \def\cal{\fam\tw@}%
  \textfont3=\tenex \scriptfont3=\tenex \scriptscriptfont3=\tenex
  \textfont\itfam=\tenit
  \def\it{\fam\itfam\tenit}%
  \textfont\slfam=\tensl
  \def\sl{\fam\slfam\tensl}%
  \textfont\bffam=\tenbf \scriptfont\bffam=\sevenbf
  \scriptscriptfont\bffam=\fivebf
  \def\bf{\fam\bffam\tenbf}%
  \textfont\ttfam=\tentt
  \def\tt{\fam\ttfam\tentt}%
  \textfont\cscfam=\tencsc
  \def\csc{\fam\cscfam\tencsc}%
  \def\big##1{{\hbox{$\left##1\vbox to8.5\p@{}\right.\n@space$}}}
  \def\Big##1{{\hbox{$\left##1\vbox to11.5\p@{}\right.\n@space$}}}
  \def\bigg##1{{\hbox{$\left##1\vbox to14.5\p@{}\right.\n@space$}}}
  \def\Bigg##1{{\hbox{$\left##1\vbox to17.5\p@{}\right.\n@space$}}}
  \setnormalsp@cing
  \usertenpointmacro
  }%
\let\userninepointmacro=\relax
\def\ninepoint{\let\pointsize=\ninepoint
  \textfont0=\ninerm \scriptfont0=\sevenrm \scriptscriptfont0=\fiverm
  \def\rm{\fam\z@\ninerm}%
  \textfont1=\ninei \scriptfont1=\seveni \scriptscriptfont1=\fivei
  \def\mit{\fam\@ne}\def\oldstyle{\fam\@ne\ninei}%
  \textfont2=\ninesy \scriptfont2=\sevensy \scriptscriptfont2=\fivesy
  \def\cal{\fam\tw@}%
  \textfont3=\tenex \scriptfont3=\tenex \scriptscriptfont3=\tenex
  \textfont\itfam=\nineit
  \def\it{\fam\itfam\nineit}%
  \textfont\slfam=\ninesl
  \def\sl{\fam\slfam\ninesl}%
  \textfont\bffam=\ninebf \scriptfont\bffam=\sevenbf
  \scriptscriptfont\bffam=\fivebf
  \def\bf{\fam\bffam\ninebf}%
  \textfont\ttfam=\ninett
  \def\tt{\fam\ttfam\ninett}%
  \textfont\cscfam=\ninecsc
  \def\csc{\fam\cscfam\ninecsc}%
  \def\big##1{{\hbox{$\left##1\vbox to8.5\p@{}\right.\n@space$}}}
  \def\Big##1{{\hbox{$\left##1\vbox to11.5\p@{}\right.\n@space$}}}
  \def\bigg##1{{\hbox{$\left##1\vbox to14.5\p@{}\right.\n@space$}}}
  \def\Bigg##1{{\hbox{$\left##1\vbox to17.5\p@{}\right.\n@space$}}}
  \setnormalsp@cing
  \userninepointmacro
  }%
\let\usereightpointmacro=\relax
\def\eightpoint{\let\pointsize=\eightpoint
  \textfont0=\eightrm \scriptfont0=\sixrm \scriptscriptfont0=\fiverm
  \def\rm{\fam\z@\eightrm}%
  \textfont1=\eighti \scriptfont1=\sixi \scriptscriptfont1=\fivei
  \def\mit{\fam\@ne}\def\oldstyle{\fam\@ne\eighti}%
  \textfont2=\eightsy \scriptfont2=\sixsy \scriptscriptfont2=\fivesy
  \def\cal{\fam\tw@}%
  \textfont3=\tenex \scriptfont3=\tenex \scriptscriptfont3=\tenex
  \textfont\itfam=\eightit
  \def\it{\fam\itfam\eightit}%
  \textfont\slfam=\eightsl
  \def\sl{\fam\slfam\eightsl}%
  \textfont\bffam=\eightbf \scriptfont\bffam=\sixbf
  \scriptscriptfont\bffam=\fivebf
  \def\bf{\fam\bffam\eightbf}%
  \textfont\ttfam=\eighttt
  \def\tt{\fam\ttfam\eighttt}%
  \textfont\cscfam=\eightcsc
  \def\csc{\fam\cscfam\eightcsc}%
  \def\big##1{{\hbox{$\left##1\vbox to8.5\p@{}\right.\n@space$}}}
  \def\Big##1{{\hbox{$\left##1\vbox to11.5\p@{}\right.\n@space$}}}
  \def\bigg##1{{\hbox{$\left##1\vbox to14.5\p@{}\right.\n@space$}}}
  \def\Bigg##1{{\hbox{$\left##1\vbox to17.5\p@{}\right.\n@space$}}}
  \setnormalsp@cing
  \usereightpointmacro
  }%
%
%
%
%
\expandafter\ifx \csname cscfam\endcsname \relax \newfam\cscfam \fi
\newtoks\baselinefactor\baselinefactor={1.2}%
\def\setnormalbaselines{%
  \normalbaselineskip=\the\baselinefactor em\relax
  }%
\def\setnormalsp@cing{
  \rm                              
  \setnormalbaselines
  \normalbaselines
  \abovedisplayskip=1.2em plus .3em minus .9em%
  \abovedisplayshortskip=0em plus .3em%
  \belowdisplayskip=1.2em plus .3em minus .9em%
  \belowdisplayshortskip=.7em plus .3em minus .4em%
  \setbox\strutbox=\hbox{\vrule height .7\baselineskip
                                depth .3\baselineskip width\z@}%
  }%
\catcode`@=12 
\edef\fmtname{\fmtname-augmented}       
\edef\fmtversion{\fmtversion-CM861707}  
%
%
%
\centerline{\hsize=12cm 
\vbox{\noindent
\ninepoint {\sl Abstract.}  In this paper we give a
characterization of the height of K3 surfaces in characteristic $p>0$.
This enables us to calculate the cycle classes of the loci in families
of K3 surfaces where the height is at least $h$. The formulas for such
loci can be seen as generalizations of the famous formula of Deuring
for the number of supersingular elliptic curves in characteristic $p$.
In order to describe the tangent spaces to these loci we study the
first cohomology of higher closed forms.} } 
\bigskip
\def\chno{1}
\centerline{\bf \chno.  Introduction}
\smallskip
\noindent
Elliptic curves in characteristic $p$ come in two sorts: ordinary and
supersingular. The distinction can be expressed in terms of the formal
group of an elliptic curve. Multiplication by $p$ on the formal group
takes the form
$$
[p](t)= at^{p^h} + \hbox{ higher order terms},\eqno(1)
$$
where $a\neq 0$ and $t$ is a local parameter. The number $h$ satisfies $1\leq h
\leq 2$ and is called the {\sl height}. By definition, the elliptic curve is
ordinary if $h=1$ and supersingular if $h=2$. There is a classical formula of
Deuring  for the number of supersingular elliptic curves over an
algebraically closed field $k$ of characteristic $p$:
$$
\sum_{\hbox{ $E$ supers.\ $/ \cong$ }}{1 \over \# {\rm Aut}(E)} =
{p-1\over 24},
$$
where the sum is over supersingular elliptic curves over
$k$ up to isomorphism. 
\par
If one views K3 surfaces as a generalization of
elliptic curves, one can make a similar distinction of K3 surfaces in
characteristic $p$ by using the formal Brauer group as Artin showed. The formal
Brauer group is a $1$-dimensional formal group associated to the second \'etale 
cohomology with coefficients in the multiplicative group. Multiplication by $p$
in this formal group has the form (1), but now we have $1\leq h \leq 10$ or
$h=\infty$, the latter if multiplication by
$p$ is zero. The height can be used to define a stratification of the
moduli spaces of K3 surfaces. A generic K3 surface will have $h=1$;
those with
$h=\infty$ are most special in this respect and called supersingular. 
\par
In this paper we first express the height of a K3 surface in terms of the action
of the Frobenius morphism on the second cohomology group with coefficients in
the sheaf $W(O_X)$ of Witt vectors of the structure sheaf $O_X$. The natural
co-filtration $W_n(O_X)$ of $W(O_X)$ induces co-filtrations on the cohomology
which correspond to approximations of the formal group. Using this
characterization we can calculate the cycle classes of the strata in the moduli
space where the height $\geq h$. This is done by interpreting the loci as
degeneracy loci of maps between bundles. The resulting formulas can be viewed as
a generalization of Deuring's formula. Generalizations of Deuring's formula to
principally polarized abelian varieties were worked out in joint work of Ekedahl
and one of us and can be found in [G]. The supersingular locus comes with a
multiplicity.
\par
In order to describe the tangent spaces to our strata we use differential forms
rather than crystalline cohomology. We calculate the dimensions of cohomology
groups $H^1(Z_i)$ and $H^1(B_i)$ where the sheaves $Z_i$ and $B_i$ are the
sheaves of certain closed
 forms introduced by Illusie. We study the dimensions of the cohomology groups
$H^1(Z_i)$ and $H^1(B_i)$ and of their images in $H^1(X,\Omega^1)$. We think
that these spaces are quite helpful to understand the geometry of surfaces in
characteristic $p$. 
\bigskip
\centerline{\hsize=12cm 
\vbox{\noindent
\ninepoint {\sl Acknowledgements.} Both authors would like to thank the
Max-Planck-Institut in Bonn for excellent working conditions during their visit
in 1998/99.  The second author would like to thank NWO and the University of
Amsterdam for their support and hospitality during his stay in Amsterdam in 1997.
We would like to thank Arthur Ogus for sharing some of his ideas with us.
} } 
\bigskip
\def\chno{1}
\centerline{\bf \chno.  Witt vector cohomology}
\smallskip
\noindent
Let $X$ be a non-singular complete variety defined over an
algebraically closed field $k$ of characteristic $p>0$. We denote by
$W_n = W_n(O_X)$ the sheaf of Witt rings of length $n$ as defined by J.-P.\
Serre, cf.\ [S]. The sheaf $W_n(O_X)$  is a coherent sheaf of rings
which comes with three operators:
{\item{i)} Frobenius $F: W_n(O_X) \to W_n(O_X)$, 
\item{ii)} Verschiebung $V: W_n(O_X) \to W_{n+1}(O_X)$,
\item{iii)} Restriction $R:W_{n+1}(O_X)\to W_n(O_X)$,}
\par
\noindent
defined by the formulas
$$
\eqalign{
F(a_0,a_1,\ldots,a_{n-1})&=(a_0^p,a_1^p,\ldots,a_{n-1}^p),\cr
V(a_0,a_1,\ldots,a_{n-1}) &= (0,a_0,a_1,\ldots, a_{n-1}),\cr
R(a_0,a_1,\ldots,a_{n})&= (a_0,a_1,\ldots,a_{n-1}).\cr
}
$$
They satisfy the relations 
$$
RVF=FRV=RFV=p.
$$
The cohomology groups $H^i(X,W_n(O_X))$ are finitely generated 
$W_n(k)$-modules. The projective system $\{ W_n(O_X), R\}_{n=1,2,\ldots}$
induces a sequence
$$
\ldots \longleftarrow H^i(X,W_n(O_X)) {\buildrel R \over \longleftarrow}
H^i(W_{n+1}(O_X))
\longleftarrow \ldots
$$
so that we can  define
$$
H^i(X,W(O_X)) = \hbox{\rm proj.\ lim } H^i(X,W_n(O_X)).
$$
This is a $W(k)$-module, but not necessarily a finitely
generated $W(k)$-module, cf. Section 3. The semi-linear operators $F$
and $V$ act on it and they satisfy the relations $FV=VF=p$.
\bigskip
\def\chno{2}
\centerline{\bf \chno. Formal Groups}
\smallskip
\noindent
Smooth formal Lie groups of dimension $1$ over an algebraically
closed field
$k$ of characteristic $\neq 0$  are characterized by their {\sl
height}, cf.\ [H], [Ma]. 
To a smooth formal Lie group $\Phi$ of dimension one
one can associate its {\sl covariant} Dieudonn\'e module $M=D(\Phi)$, a
free
$W(k)$-module.  It possesses two operators $F$ and $V$ with the
following properties: the operator
$F$ is
$\sigma$-linear, the operator  $V$ is $\sigma^{-1}$-linear
 and topologically nilpotent and  they satisfy $FV=VF=p$. Here $\sigma$ denotes the Frobenius map on $k$.  Then
$M$ is a free $W(k)$-module with the following properties:
\par
\item{a)}$\quad\dim (\Phi) = \dim_k(M/VM)$,
\par
\item{b)}$\quad {\rm height}(\Phi)  = {\rm rank}_W(M).$
\par\noindent
Note that one has the equalities
$$ 
{\rm rank}_W(M)=\dim_k(M/pM)
=\dim_k(M/FM)+\dim_k(M/VM).
$$
\bigskip
\def\chno{3}\centerline{\bf \chno. The Formal Brauer Group of Artin-Mazur}
\smallskip
\noindent
For a proper variety $X/k$ one may consider the formal completion of
the Picard group. The group of $S$-valued points of 
$\widehat{\rm Pic}(X)$ fits into the exact sequence
$$
0 \tto \widehat{\rm Pic}(S) \tto H^1(X\times S, \G_m) \tto
H^1(X,\G_m)
$$
for any  local artinian scheme $S$ with residue field $k$. Here
cohomology is \'etale cohomology.  This idea of studying
infinitesimal properties of cohomology was generalized to the higher
cohomology groups $H^r(X,\G_m)$ by Artin and Mazur, cf.\ [A-M]. Their work
leads to contravariant functors $\Phi^r\colon {\rm Art} \to {\rm Ab}$
with
$$
\Phi^r(S)=\ker H^r(X\times S, \G_m) \tto H^r(X,\G_m),
$$
which under suitable circumstances are representable by formal Lie
groups. For a K3 surface $X$ this is the case and we find for $r=2$ the
formal Brauer group $\Phi = \Phi_X = \Phi^2$. Its tangent space is
$$
T_{\Phi}= H^2(X,O_X).
$$
For a K3-surface $X$ we have two possibilities:
\par
\indent i) $h(\Phi)=\infty$ and $\Phi = \hat{\G}_a$,  the formal
additive group. The K3-surface is called {\sl supersingular} (in the
sense of Artin).
\hfill\break \indent ii) $h(\Phi)< \infty$. Then $\Phi$ is a
$p$-divisible formal group. Moreover, it is known that $1 \leq h(\Phi)
\leq 10$. This follows from the following theorem of Artin, cf.\
[A]. We shall write simply $h$ for $h(\Phi)$.
\par
\proclaim (\chno.1) Theorem. If the formal Brauer group $\Phi_X $  
of a K3 surface $X$ is $p$-divisible then it satisfies the relation
$2h \leq B_2-\rho$, where $B_2$  is the second Betti number and $\rho$ the rank
of the N\'eron-Severi group.
\par
\noindent
For the proof one combines Theorem (0.1) of [A] with Deligne's [D] result on
lifting K3 surfaces, see also [I]. We give a proof in Section 10.
 This theorem implies that if $\rho = 22$  
then necessarily we have $h=\infty$. If $h\neq \infty$ then it follows that 
$1\leq h \leq 10$. One should view $h=1$ as the
generic case.
It was conjectured by Artin that if $h=\infty$ then $\rho=22$. This is
known for elliptic K3 surfaces, see [A]. Note that a surface with 
$\rho=22$ is called {\sl supersingular} by Shioda, cf.\ [Sh].
\par
The following result by Artin and Mazur is crucial:
\proclaim (\chno.2) Theorem.  The Dieudonn\'e module of the formal Brauer
group $\Phi_X$ is given by
$$
D(\Phi_X) \cong H^2(X,W(O_X)).
$$
\par
\noindent
For the proof we refer to [A-M]. The point to notice is that 
$$
D(\Phi_X)= H^2(X,D\G_m)=H^2(X,W(O_X)).
$$ 
\noindent
{\bf (\chno.3) Remark.} Note that this explains why the Witt vector
cohomology is sometimes not finitely generated:  if $\Phi_X \cong
{\hat \G}_a$ then $H^2(X, W(O_X))$ is not finitely generated over
$W(k)$ because
$D(\hat{\G}_a)= W(k)[[T]]$.
\bigskip
\def\chno{4}
\centerline{\bf \chno. Vanishing of Cohomology}
\smallskip
\noindent
We collect a number of results on the vanishing of cohomology groups
for K3 surfaces that we need in the sequel.
\smallskip
\proclaim (\chno.1) Lemma. Let
$X$ be a K3 surface. We have
$H^1(X,W_n(O_X))=0$ for all $n>0$, hence $H^1(X,W(O_X))=0$.
\par
\noindent
{\sl Proof.} Since $X$ is a K3 surface we have by definition 
$H^1(X,O_X)=0$. The lemma is deduced from this by induction on $n$.
Assume that $H^1(X, W_{n-1}(O_X))=0$.  Then the exact sequence 
$$
0 \tto W_{n-1}(O_X){\buildrel V \over \tto} W_n(O_X){\buildrel R^{n-1}
\over
\tto} O_X\tto 0
$$
induces an exact sequence
$$
H^1(X,W_{n-1}(O_X)){\buildrel V \over \tto} H^1(W_n(O_X)) \mapright{R^{n-1}} 
H^1(X,O_X).
$$
This implies that $H^1(X,W_n(O_X))=0$. $\square$
\par
\proclaim (\chno.2) Lemma. For a projective surface $X$ with $H^1(X,O_X)=0$ the
induced map  $R\colon H^2(X,W_n(O_X))\to H^2(X,W_{n-1}(O_X))$ is surjective with
kernel $\cong H^2(X,O_X)$.
\par
\noindent
{\sl Proof.} This follows from the exact sequence
$$
0\to O_X \tto W_n(O_X) {\buildrel R \over \tto} W_{n-1}(O_X) \to 0
$$
and the vanishing of $H^1(X,O_X)$ and of $H^3(X,O_X)$. $\square$
\par
\proclaim (\chno.3) Lemma. In $H^2(X,W_n(O_X))$ we have
$$
RV(H^2(X,W_n(O_X))) = V(H^2(X,W_{n-1}(O_X))).
$$ 
\par
\noindent
{\sl Proof.} The commutativity of the diagram
$$
\matrix{ W_n(O_X) & \mapright{V}& W_{n+1}(O_X)\cr
\mapdown{R} && \mapdown{R} \cr
W_{n-1}(O_X) & \mapright{V} & W_n(O_X).\cr
}
$$
gives in cohomology a commutative diagram
$$
\matrix{ H^2(X,W_n(O_X)) & \mapright{V}& H^2(X,W_{n+1}(O_X))\cr
\mapdown{R} && \mapdown{R} \cr
H^2(X,W_{n-1}(O_X)) & \mapright{V} & H^2(X,W_n(O_X)).\cr
}
$$
The surjectivity of the left hand $R$, which follows from the
preceding lemma, implies the claim. 
$\square$
\proclaim (\chno.4) Lemma. Assume that for some  $n>0$ the  map
$F: H^2(X,W_n(O_X)) \tto H^2(X,W_n(O_X))$ vanishes. Then for all $0\leq
i\leq n$ the map $F: H^2(X,W_i(O_X)) \tto H^2(X,W_i(O_X))$ is zero.
Moreover, for all $0\leq i \leq n$ the module $H^2(X,W_i(O_X))$ is a
vector space over $k$.\par
\noindent
{\sl Proof.} The first result follows from the commutativity of the
diagram
$$
\matrix{ H^2(X,W_n(O_X)) & \mapright{R^{n-i}}& H^2(X,W_{i}(O_X))\cr
\mapdown{F} && \mapdown{F} \cr
H^2(X,W_{n}(O_X)) & \mapright{R^{n-i}} & H^2(X,W_i(O_X))\cr
}
$$
and Lemma (\chno.2). The second claim follows from $p=FVR$ and $k\cong
W_i(k)/pW_i(k)$. $\square$
\proclaim (\chno.5) Lemma. Assume that $X$ is a K3 surface. 
The following two sequences are exact:
$$
 0 \to H^2(X,W_{n-1}(O_X)) \mapright{V} H^2(X,W_n(O_X))
\mapright{R^{n-1}} H^2(X,O_X) \to 0,
$$
$$
0 \to H^2(X,W(O_X)) \mapright{V} H^2(X,W(O_X))
\mapright{R'} H^2(X,O_X) \to 0,
$$
where $R'$ is the map induced by $W_n(O_X) {\buildrel R^{n-1} \over 
\longrightarrow} W_1(O_X)$ as $n \to \infty$.
\par
\noindent
{\sl Proof.} The first exact sequence follows from the exact sequence
$$
0 \to W_{n-1}(O_X) \mapright{V} W_n(O_X) \mapright{R^{n-1}} O_X \to 0
$$
and Lemma (4.2). Because the projective system $H^2(X,W_n(O_X))$
satisfies the Mittag-Leffler condition we may take the projective
limit.
$\square$
\bigskip
\def\chno{5}
\centerline{\bf \chno. Characterization of the Height}
\smallskip
\noindent
Let $X$ be a K3 surface and let $\Phi_X$ be its formal Brauer group in
the sense of Artin-Mazur. The isomorphism class of this formal group
is determined by its height $h$. The following theorem expresses this
height in terms of Witt vector cohomology.
\proclaim (\chno.1) Theorem. The height satisfies  $h(\Phi_X) \geq i+1$  if and
only if the Frobenius map $ F\colon H^2(X,W_i(O_X)) \to H^2(X,W_i(O_X))$ is
the zero map.
\par
\proclaim (\chno.2) Corollary. We have the following characterization of the
height:
$$
h(\Phi_X) = \min \{ i \geq 1 \colon [F: H^2(W_i(O_X)) \to
H^2(W_i(O_X))] \neq 0 \}.
$$\par
\smallskip
\noindent
{\sl Proof of the Theorem.} ``$\Leftarrow$"  In case $
h(\Phi_X)= \infty$ the implication $\Leftarrow$ is immediate. So
we may consider the case where the height of $\Phi_X$ is finite. 
Assume that the map $F\colon H^2(X,W_i(O_X)) \to H^2(X,W_i(O_X))$ is the zero
map. We set
$$
M=D(\Phi)\cong H^2(X,W(O_X)), \quad {\hbox { \rm the covariant
Dieudonn\'e module.}}
$$
Since $\dim_k(H^2(X,W(O_X))/VH^2(X,W(O_X)) =1$ by Lemma (4.5),  we have
by b) in Section 2
$$
\dim_k(H^2(X,W(O_X))/FH^2(X,W(O_X))=h-1.
$$
The surjectivity of the projection $H^2(X,W(O_X)) \tto H^2(X,W_i(O_X))$
implies the surjectivity of
$$
H^2(X,W(O_X))/FH^2(X,W(O_X)) \tto H^2(X,W_i(O_X))/FH^2(X,W_i(O_X)).
$$
By assumption we have
$H^2(X,W_i(O_X))/FH^2(X,W_i(O_X))\cong H^2(X,W_i(O_X))$ and by Lemma
(4.5) we have
$$
\dim_kH^2(X,W_i(O_X))=i,
$$
i.e.\ we find $h-1 \geq i$, or equivalently, $h\geq i+1$.
\par
Conversely, we now  prove ``$\Rightarrow$". If $h(\Phi_X) =
\infty$ then $\Phi_X = \hat \G_a$, the formal additive group of
dimension $1$. So $F$ acts as zero on $D(\hat
\G_a)=D(\Phi_X)=H^2(X,W(O_X))$. As in Lemma (4.4) we conclude that $F$
acts on $H^2(X,W_i(O_X))$ as the zero map. Therefore we may assume
that $h(\Phi_X)=h < \infty$. We thus assume that $h(\Phi_X) \geq i+1$. 
We set
$$
H= H^2(X,W(O_X))
$$
and have 
$$
V^{h-1}H \subset \ldots \subset V^2H \subset VH \subset H.
$$
Under projection this is mapped surjectively to
$$
0 \subset V^{h-2}H^2(O_X) \subset \ldots \subset VH^2(W_{h-2}(O_X)) \subset
H^2(X,W_{h-1}(O_X)).
$$
All the inclusions are strict because of Lemma (4.5).
\proclaim Claim. We have  $V^{h-1}H^2(X,W(O_X)) = FH^2(X,W(O_X))$.\par
\noindent
{\sl Proof of the claim.}  Since our modules are free over $W$ we deduce  from
Manin's results~[M] (but see also [H] because we use the covariant theory):
$$
D(\Phi_X)\cong W[F,V]/W[F,V](F-V^{h-1}).
$$
Note that $F-V^{h-1}$ is written on the right. But we can transfer it to the
left  using $FV=p=VF$ as follows:
$$
\eqalign{F(\sum a_{ij}F^iV^j)&=\sum a_{ij}^{\sigma} FF^iV^j=\sum
(a_{ij}^{\sigma} F^iV^j)F=\cr &=\sum (a_{ij}^{\sigma}F^iV^j)V^{h-1} =
V^{h-1}(\sum a_{ij}^{\sigma^h} F^iV^j).\cr}
$$
This together with Theorem (3.2) proves the claim.
\par
We now find $FH^2(W_{h-1}(O_X))=0$. By Lemma (4.4) we conclude that $F$ acts
on
$H^2(W_i(O_X))$ for $i \leq h-1$ as zero.
$\square$
\par
\proclaim (\chno.3) Corollary. The height of $\Phi_X$ is $\infty$
if and only if  the Frobenius endomorphism  
$F: H^2(X,W_{10}(O_X)) \to H^2(X,W_{10}(O_X))$ is zero. 
\par
\noindent
{\sl Proof.} If the height is finite, then we know by Artin and Mazur (see
(3.1)) that we have
$h\leq 10$. $\square$

\proclaim
(\chno.4) Corollary. Set $H=H^2(X,W_{10}(O_X))$ and consider the filtration
$$
\{ 0 \} \subset R^9V^9H \subset R^8V^8H \subset \ldots \subset
R^{h-1}V^{h-1}H \subset \ldots \subset H.$$
If $h$ is the height of $\Phi_X$ then $F(H)= R^{h-1}V^{h-1}(H)$.
\par
\noindent
{\sl Proof.} The $(h-1)$-th step $V^{h-1}H^2(W(O_X))$ in the filtration
$$
V^{10}H^2(W(O_X)) \subset V^9H^2(W(O_X)) \subset \ldots \subset H^2(W(O_X))
$$
maps surjectively to the corresponding step $R^{h-1}V^{h-1}H$ of the
filtration on $H$. By our claim we have
$$
V^{h-1}H^2(W(O_X)) = FH^2(W(O_X)).
$$ 
This implies the assertion. $\square$
\proclaim (\chno .5) Corollary. If $h(\Phi_X)=h< \infty$ and if
$\{ \omega,V\omega,V^2\omega, \ldots, V^{h-1}\omega\} $ is a $W$-basis of
$H^2(X,W(O_X))$ then $F$ acts as zero on $H^2(X,W_i(O_X))$ if and only if
$F(\bar\omega)=0$, with $\bar\omega$ the image of $\omega$ in
$H^2(X,W_i(O_X))$.\par
\proclaim (\chno .6) Corollary. If $h(\Phi_X)=h<\infty$,
then $\dim_{k} \ker[F: H^2(W_i) \to H^2(W_i)] = \min \{i,h-1\}$. \par
\noindent
{\sl Proof.} By Lemma (4.5) and Corollary (5.2), we have 
$\dim_{k} \ker[F: H^2(W_i) \to H^2(W_i)]$ $= i$ if 
$i \leq h -1$.
Assume $i \geq h$. Using the notation in Corollary 5.5, we know that 
$\langle V^{i-h +1}R^{i-h +1}{\bar\omega},V^{i-h +2}R^{i-h +2}{\bar\omega},
V^{i-h +3}R^{i-h +3}{\bar\omega}, \ldots, V^{i-1}R^{i-1}{\bar\omega}\rangle$ 
is a basis of $\ker F$.
$\square$
\smallskip
The case $h\geq 2$ is characterized by the vanishing of Frobenius on $H^2(O_X)$.
We now formulate in an inductive way a similar characterization of the condition
$h \geq n+1$. If for a K3 surface $X$ one assumes that $F$ is zero on 
$H^2(W_i)~(i = 1,  \ldots n-1)$ then one has $F H^2(W_n) \subset V^{n-1}H^2(O_X)$
and $F$ vanishes on $VH^{2}(W_{n-1})$. Since we have a natural
$({\sigma}^{-1})^{n-1}$-isomorphism $ H^2(O_X) \cong V^{n-1}H^2(O_X)$,
 one has an induced homomorphism 
$$\phi_{n} \colon H^{2}(O_{X}) \cong H^2(W_n)/VH^2(W_{n-1})
\to V^{n-1}H^2(O_X) \cong  H^2(O_X). \eqno(2)
$$
This map is ${\sigma}^{n}$-linear.
The following theorem is clear by the construction of $\phi_n$.
\proclaim (\chno.7) Theorem. Suppose $F$ is zero on $H^2(W_i)$ for 
$i = 1, \ldots n-1$. Then $F$ vanishes on $H^2(W_n)$ if and only if
$\phi_{n}: H^2(O_X) \to H^2(O_X)$ vanishes. \par
\noindent
\bigskip
\noindent
\def\chno{6}\centerline{\bf \chno. Closed Differential Forms }
\smallskip
\noindent
Let $F: X \to X^{(p)}$ be the relative Frobenius morphism of a  K3
surface $X$. By means of the Cartier operator $C: \Omega_{X,\rm
closed}^{\bullet}\to \Omega^{\bullet}_X$  we can define sheaves $B_i\Omega_X^1$
of rings inductively by $B_0\Omega_X^1=0$, $B_1\Omega_X^1=dO_X$ and
$C^{-1}(B_i\Omega_X^1) = B_{i+1}\Omega_X^1$. Similarly, we define sheaves
$Z_i\Omega_X^1$  inductively by $Z_0\Omega_X^1=\Omega_X^1$, $Z_1\Omega_X^1=
\Omega^1_{X,\rm closed}$, the sheaf of $d$-closed forms and by setting 
$$
Z_{i+1}\Omega_X^1\colon = C^{-1}(Z_i\Omega_X^1).
$$
Usually we simply write $B_i$  and $Z_i$. The sheaves $B_i$ and $Z_i$ can be 
viewed as locally free subsheaves of
$(F^i)_*\Omega_X^1$ on $X^{(p^i)}$. They were introduced by Illusie in [Il] and
can be used to  provide de Rham-cohomology with a rich structure.   The inverse
Cartier operator gives rise to an isomorphism 
$$ C^{-i} : \Omega_{X^{(p^i)}}^1 {\buildrel \simeq \over \tto}
Z_i/B_i
$$
or a $\sigma^{-i}$-linear isomorphism $\Omega_X^1 \cong Z_i/B_i$.  Note that we
have the inclusions
$$
0=B_0 \subset B_1 \subset \ldots \subset B_i \subset \ldots \subset Z_i
\subset \ldots \subset Z_{1} \subset Z_0 = \Omega_X^1.
$$
We also have an exact sequence   
$$ 
0 \to Z_{i+1} \tto Z_i {\buildrel dC^i
\over \ttto } d\Omega_X^1 \to 0. \eqno(3)
$$
\proclaim (\chno.1) Lemma. If $X$ is a K3 surface $X$ 
we have  {\rm i)} $H^{0}(B_{i}) = 0$ for all $i \geq 0$;
 {\rm ii)}  the natural inclusion $B_{i} \to B_{i + 1}$ induces an injective 
homomorphism $H^{1}(B_{i}) \to H^{1}(B_{i + 1})$. \par
\noindent
{\sl Proof.} i) The  natural injection $B_{i} \to \Omega^{1}_{X}$
induces an injection $H^0(B_i) \to H^0(\Omega^{1}_{X})$ and
we know $H^{0}(\Omega^{1}_{X}) = 0$. ii) This follows from i) and the
exact sequence
$$
 0 \to B_{i} \tto B_{i+1}  {\buildrel C^i \over \tto} B_{1} \to 0
$$
with $C$ the Cartier operator. $\square$  
\par
\bigskip
There is a close relationship between the Witt vector cohomology and the
cohomology of $B_i$ as follows. Serre introduced in [S] a map
$D_i: W_i(O_X)
\tto
\Omega_X^1$ of sheaves in the following way: 
$$
D_i(a_0,a_1,\ldots,a_{i-1}) =  a_0^{p^{i-1}-1}da_0 +\ldots + 
a_{i-2}^{p-1}da_{i-2} +
da_{i-1}.
$$
It satisfies $D_{i+1}V = D_{i}$,
and Serre showed that this induces an injective map of sheaves of additive groups
$$
D_i : W_i(O_X)/FW_i(O_X) \tto \Omega^1_X 
$$
inducing an isomorphism 
$$
D_i: W_i(O_X)/FW_i(O_X) {\buildrel \sim \over \tto} B_i\Omega_X^1 \eqno(4)
$$
The exact sequence  $0 \to W_i {\buildrel F \over \tto} 
W_i \tto W_i/FW_i \to 0$ gives rise to the  exact sequence
$$
0 \to H^1(W_i/FW_i) \to H^2(W_i) {\buildrel F \over \tto}
H^2(W_i) \to H^2(W_i/FW_i) \to 0 \eqno(5)
$$
and we thus  have an isomorphism $H^1(W_i/FW_i) \cong
\ker[F: H^2(W_i) \to H^2(W_i)]$.
Combining the result on the dimension of the kernel of $F$ on $H^2(W_i)$  from
Section 5 with  (4) we get an interpretation of the height $h$ in terms of
the groups $H^1(B_i)$. 
\par
\proclaim (\chno.2) Theorem. We have  
$$
 \dim H^1(B_i)  = \cases{ \min \{i,h-1\} &if $h \neq \infty$,\cr
i & if $h=\infty$. \cr}
$$
\par
The Verschiebung induces an exact sequence
$$
0 \to W_i/FW_i {\buildrel V\over \tto} W_{i+1}/FW_{i+1} \to
O_X/FO_X \to 0 
$$
and this gives rise to
$$
0 \to H^1(W_i/FW_i) {\buildrel V \over \tto } H^1(W_{i+1}/FW_{i+1}) \tto
H^1(O_X /FO_X )  \to \ldots 
$$ 
i.e., Verschiebung induces for all $i$ an {\sl injective map}. Moreover,  it is
surjective if and only if $h \neq \infty$ and $i\geq h-1$.
\par
We have a commutative diagram (with $\beta_i$ the natural map induced by
 $B_i \subset \Omega^{1}_{X}$)
$$
\matrix{ H^1(W_i/FW_i) & {\buildrel D_i \over \tto} & H^1(\Omega^{1}_{X}) \cr
\mapdown{\cong} && \mapdown{=}\cr
H^1(B_i\Omega^{1}_{X}) & {\buildrel {\beta_i} \over \ttto } & H^1(\Omega^{1}_{X}) \cr}
$$
We study  the kernel of $D_i$, equivalently  the
kernel of the natural map $\beta_i: H^1(B_i\Omega^{1}_{X}) 
\to H^1(\Omega^{1}_{X})$ in Sections 9-11. 
\par
\smallskip
\noindent
\proclaim (\chno.3) Lemma. The Euler-Poincar\'e characteristics of $B_i$ and
$Z_i$ are given by $\chi(B_i)=0$ and $\chi(Z_i)=-20$.\par 
\noindent {\sl Proof.}  Since the kernel and the cokernel of $F$ on $H^2(W_i)$
have the same dimension by (5) the result for $B_i$ follows from (4) and (5). The
identity
$\chi(B_i) + \chi(\Omega^{1}_{X})=\chi(Z_i)$ resulting from the isomorphism 
$Z_i/B_i \cong \Omega^{1}_{X}$implies the result. $\square$
\par
\bigskip
\noindent
\def\chno{7}\centerline{\bf \chno. De Rham Cohomology}
\smallskip
\noindent
The de Rham cohomology of a K3 surface is the hypercohomology of the
complex $(\Omega_X^{\bullet},d)$. The dimensions $h^{p,q}$ of the graded
pieces are given by the Hodge diamond.
$$
\matrix{ && 1 && \cr
&0&&0& \cr
1&&20 && 1\cr
&0&&0&\cr
&&1 &&\cr}
$$
On $H^2_{dR}$ we have a perfect pairing $\langle \, , \, \rangle$ given by
Poincar\'e duality; cf.\ [D].
\par
The Hodge spectral sequence with $E_1^{ij}=H^j(X, \Omega^{i}_{X})$ converges
to $H_{dR}^*(X)$. The second spectral sequence of hypercohomology has
as $E_2$-term $E_2^{ij}=H^i({\cal H}^j(\Omega^{\bullet}_{X}))$ abutting to
$H_{\rm dR}^{i+j}(X/k)$. But the Cartier operator yields an isomorphism
of sheaves  
$$
C^{-1}: \Omega^i_{X^{(p)}} {\buildrel \sim \over \tto} {\cal
H}^i(F_*(\Omega_{X/k}^{\bullet})),
$$
so that we can rewrite this as
$$
E_2^{ij}=H^i(X^{\prime}, {\cal H}^j(\Omega^{\bullet}))\cong  H^i(X^{\prime},
\Omega_{X^{\prime}}^j) \Rightarrow H_{dR}^*(X),
$$
where $X^{\prime}=X^{(p)}$ is the base change of $X$ under Frobenius. We thus
get two filtrations on the de Rham cohomology: the Hodge filtration
$$
(0) \subset F^2 \subset F^1 \subset H_{dR}^2,
$$
and the so-called {\sl conjugate filtration}
$$
(0) \subset G_1 \subset G_2 \subset H_{dR}^2.
$$
We have $\rk(F^1)=\rk(G_2)=21$, $\rk(F^2)=\rk(G_1)=1$ and
$$
(F^1)^{\bot}=F^2\qquad{\rm and}\qquad G_1^{\bot}=G_2.
$$
We have also
$$
F^1/F^2 \cong H^{1}(X, \Omega^{1}_{X}), \quad G_2/G_1 \cong
H^1(X,\Omega^{1}_X),
$$
cf.\ [D]. Moreover, from the description with the second spectral sequence it
follows that $G_1$ is the image under Frobenius of $H^2_{\rm dR}$ and also of
$H^2(X,O_X)$. The conjugate filtration is an analogue of the complex conjugate
of the Hodge filtration in characteristic zero. 
\par
The relative position of these two filtrations is an interesting
invariant of a K3 surface. We have the three cases
\item{a)} $F^1\cap G_1= \{ 0\}$;\par
\item{b)} $G_1 \subset F^1$; $G_1 \neq F^2$;\par
\item{c)} $G_1=F^2$.\par
\noindent
The first case happens if and only if $F: H^2(X,O_X) \to H^2_{\rm dR}(X) 
\to H^2(X,O_X)$ is not zero, i.e.\ if $h=1$. Such $X$ are called {\sl ordinary}.
The second case happens if $h\geq 2$, while the last case is by definition
the {\sl superspecial} case. In this case the two filtrations coincide. 
It is known that two superspecial K3
surfaces are isomorphic (as unpolarized varieties)(cf.\ [O]).
\par
We have the following result of Ogus (cf.\ [O]) which provides us with an
interpretation of $H^1(Z_1)$.
\par
\proclaim (\chno.1) Proposition. We have an isomorphism $
F^1 \cap G_{2} \cong  H^{1}(X, Z_1)$. 
\par
\noindent
{\sl Proof.} The map $z_1: H^1(Z_1) \to H^2_{\rm dR}$ given by
$\{f_{ij}\}\mapsto (0,\{f_{ij}\},0)$ is injective. Indeed, if  $\{f_{ij}\}$
represents an element in the kernel, then it is of the form
$(\delta h_{ij}, dh_{ij} + \omega_j-\omega_i, d\omega_i)$
for a $h_{ij} \in C_1(O_{X})$, $\omega_i \in C_0(\Omega^{1}_{X})$. Then the
$\omega_i$ are closed and $h_{ij}$ defines a cocycle. Since $H^1(O_{X})=0$
we can write $dh_{ij}=\eta_j-\eta_i$ and $f_{ij}$ is a coboundary.
The image is contained in $F^1$ and is orthogonal to the
image $G_1$ of Frobenius. Indeed, take a class $F(a)$ and consider the 
cupproduct $\langle F(a), z_1(f_{ij})\rangle$. Applying the Cartier operator
we see that it is zero. But $C: H^4_{\rm dR} \to H^4_{\rm dR}$ is
a bijection. Hence the image lies in $F^1 \cap G_2$. This implies that
$\dim H^1(Z_1) \leq 20$ if $X$ is not superspecial and $\leq 21$ for 
superspecial $X$. The exact sequence
$$
0\to H^1(B_1) \to H^1(Z_1) {\buildrel C \over \tto} H^1(\Omega^{1}_{X})
 \to H^2(B_1) \to H^2(Z_1)\to 0
$$
implies together with the value of $h^1(B_1)=h^2(B_1)$  and $\chi(Z_1)=-20$
that $h^1(Z_1)=20$ unless $X$ is superspecial. 
But if $X$ is superspecial then because of $F^2=G_1$ the Cartier
operator gives an isomorphism $C: H^1(Z_1)/H^1(B_1) \cong
H^1(\Omega^{1}_{X})$ implying that $h^1(Z_1)=21$. $\square$
\bigskip
\noindent
\def\chno{8}\centerline{\bf \chno. An Extension of de Rham Cohomology}
\smallskip
\noindent
We define an extension of de Rham cohomology by considering an enlarged
complex. (It captures the $[0,1)$-part of crystalline cohomology.) It is
defined as follows. We denote by
$H^{i}_{\rm dRW_{n}}(X/S)$ the cohomology  of the double complex $CW_{n}$ of
additive groups which is defined by the commutative digram:
$$
\def\normalbaselines{\baselineskip20pt
\lineskip3pt \lineskiplimit3pt }
\matrix
{{\partial}\uparrow  &  & {\partial}\uparrow &  & 
{\partial}\uparrow \cr
C_{2}(W_{n}(O_{X/S})) & {\buildrel D_{n} \over \tto} &
C_{2}(\Omega_{X/S}^{1}) & {\buildrel d \over \tto}
& C_{2}(\Omega_{X/S}^{2})  \cr
{\partial}\uparrow  &  & {\partial}\uparrow &  & 
{\partial}\uparrow \cr
C_{1}(W_{n}(O_{X/S})) & {\buildrel D_{n} \over \tto} &
C_{1}(\Omega_{X/S}^{1}) & {\buildrel d \over \tto}
& C_{1}(\Omega_{X/S}^{2})  \cr
{\partial}\uparrow  &  & {\partial}\uparrow &  & 
{\partial}\uparrow \cr
C_{0}(W_{n}(O_{X/S})) & {\buildrel D_{n} \over \tto} &
C_{0}(\Omega_{X/S}^{1}) & {\buildrel d \over \tto}
& C_{0}(\Omega_{X/S}^{2}),  \cr
}
$$
where $C_{i}$ are the $i$-th {\v C}ech 
cochains, $D_{n}$ are the maps induced by the differential
of Serre given in Section 6, the differentials $d$ are defined
by the exterior differentiation of differential forms and the vertical
differentials are taken in the {\v C}ech sense.
As usual, we denote by $\delta$ the differential of the
single complex associated with $CW_{n}$. 
An element of  $H^{2}_{\rm dRW_{n}}(X/S)$ is represented by 
a triple $(\alpha_{0}, \alpha_{1}, \alpha_{2}) \in  C_{2}(W_{n}(O_{X/S})) \oplus C_{1}(\Omega_{X/S}^{1})
\oplus C_{0}(\Omega_{X/S}^{2})$. In case $n = 1$, 
$H^{2}_{\rm dRW_{1}}(X/S)$ is nothing but the de Rham cohomology
$H^{2}_{\rm dR}(X/S)$. On $H^{2}_{\rm dRW_{n}}(X/S)$ we have the Hodge
filtration
$$
   0 \subset F^{2} \subset F^{1} \subset H^{2}_{\rm dRW_{n}}(X/S).
$$
Here the  $F^{i}$ (for $i\neq 0$) is naturally isomorphic to the
$F^i$-part in  the Hodge filtration of $H^{2}_{\rm dR}(X/S)$. We
have a natural isomorphism
$$
     H^{2}_{\rm dRW_{n}}(X/S)/F^{1} \cong H^{2}(X, W_{n}(O_{X/S})).
$$ 
Since the Frobenius morphism $F$ is a zero map on $F^{1}$,
we have an induced homomorphism
$$
       F : H^{2}(X, W_{n}(O_{X/S})) \tto  H^{2}_{\rm dRW_{n}}(X/S).
$$
The map
$V^{n-1} : O_{X/S} \tto W_{n}(O_{X/S})$ gives rise to
a homomorphism 
$$
V^{n - 1} : C_{i}(O_{X/S}) \tto C_{i}(W_{n}(O_{X/S})).
$$
Using this homomorphism and taking the identity mapping from 
$C_{i}(\Omega_{X/S}^{j})$ to $C_{i}(\Omega_{X/S}^{j})$, we have 
a homomorphism of complexes of additive groups $CW_{1} \tto CW_{n}$.
Therefore, we have a homomorphism of additive groups:
$$
 V^{n-1}:  H^{2}_{\rm dR}(X/S) \tto H^{2}_{\rm dRW_{n}}(X/S).
$$
Let $X_0$ be a K3 surface over a field $k$ and assume now that $F :
H^{2}(W_{i}(O_{X_{0}}))
\tto H^{2}(W_{i}(O_{X_{0}}))$  is zero for 
$i = 1, \ldots, n - 1$. Then, by the same argument as in Section~5,
we have $FH^{2}_{\rm dRW_{n}}(X_{0}) \subset V^{n-1}H^{2}_{\rm dR}(X_{0})$.  
Therefore, using the inverse of the natural isomorphism of additive groups 
$H^{2}_{\rm dR}(X_{0}) \cong V^{n-1}H^{2}_{\rm dR}(X_{0})
\subset H^2_{\rm dRW_n}(X_0)$, we have a homomorphism 
$$
  \Phi_{n} : H^{2}(W_{n}({O}_{X}))
{\buildrel F\over\longrightarrow} 
V^{n-1}H^{2}_{DR}(X)
\cong H^{2}_{DR}(X).
$$
Since we have $H^{2}(W_{n}(O_{X}))/VH^{2}(W_{n-1}(O_{X})) \cong
H^{2}(O_{X_0}) \cong H^{2}_{DR}(X)/F^{1}$, and
since $\Phi_{n}$ maps $VH^{2}(W_{n-1}(O_{X_0}))$ to $F^{1}$, the map
$\Phi_{n}$ induces a homomorphism from $H^{2}(O_{X_0})$ to $H^{2}(O_{X_0})$.
This homomorphism coincides with the map $\phi_{n}$ which was
constructed in Section 5.
\par
We now take a basis $\omega_{0}$ of
$H^{0}(X_{0},\Omega_{X_{0}}^{2})$ and take the dual basis $\zeta_{0}$ of
$H^{2}(X_{0}, O_{X_{0}})$. Via the Hodge filtration of
$H^{2}_{\rm dR}(X_{0})$, we can naturally regard
$H^{0}(X_{0},\Omega_{X_{0}}^{2})$ as a subspace of $H^{2}_{\rm dR}(X_{0})$.
Therefore, we may assume $\omega_{0}$ is an element of $H^{2}_{\rm
dR}(X_{0})$.
\par
Since $R^{n-1} : H^{2}(W_{n}(O_{X_{0}})) \rightarrow H^{2}(O_{X_0})$
is surjective, there exists an element $\alpha_{0} \in H^{2}(W_{n}(O_{X_{0}}))$
such that $R^{n-1}(\alpha_{0}) = \zeta_{0}$. Then, by Theorems 5.1 and 5.7
we have the following proposition.
\par
\proclaim
(\chno.2) Proposition. Suppose that for a K3 surface $X_{0}$ the map 
$F : H^{2}(W_{i}(O_{X_{0}})) \to H^{2}(W_{i}(O_{X_{0}}))$ 
is zero for  $i = 1, \ldots, n - 1$. Then, with the notation introduced above,
$h(\Phi_{X_{0}}) \geq n + 1$ if and only if $\langle \Phi_{n}({\alpha}_{0}),
\omega_{0} \rangle = 0$. \par
\bigskip
\def\chno{9}
\centerline{\bf \chno. The Dimensions of the Spaces of Closed Forms}
\smallskip
\noindent
We study the dimensions of the spaces $H^1(X,B_n)$ and $H^1(X,Z_n)$. We also
consider their images in $H^2_{\rm dR}(X)$ and this gives a finer structure
on these de Rham cohomology groups.
\par
Let us consider the natural map $\beta_n$ induced in $H^1$ by the inclusion $B_n
\subset \Omega^{1}_{X}$:
$$
\beta_n: H^1(B_n) \tto H^1(\Omega^1).
$$
\par
\proclaim (\chno.1) Proposition. If ${\beta}_n$ is not injective
then ${\beta}_m$ is not injective for every $m\geq n$ and  $\dim H^1(B_m) <
\dim H^1(B_{m+1})$.
\par
\noindent
{\sl Proof.} The maps $\beta_n$ are compatible with the natural maps $H^1(B_n)
\to H^1(B_{n+1})$  and by (6.1) these maps $H^1(B_n) \to H^1(B_{n+1})$ are
injective. If $\beta_n$ is not injective it follows that $\beta_{n+1}$ is not
injective. To prove the second statement, we start with the case $n=1$. If
 $\beta_1: H^1(B_1) \to H^1(\Omega^{1}_{X})$ is not injective then there
exists a non-trivial  cocycle $f_{ij} \in C^1(O_X/FO_X)$ and a 1-cochain
${\omega_i}$ of 1-forms such that  $ df_{ij}= \omega_j -\omega_i$.
Since for affine open sets $U$ the Cartier map  $H^0(\Omega_{U, \rm closed}^1)
\to H^0(\Omega^{1}_{U})$ is surjective we can find closed forms
$\tilde{\omega}_i$ and regular functions $g_{ij}$ on $U_{i} \cap U_{j}$ such
that we have a  relation
$$
f_{ij}^{p-1} df_{ij}+ dg_{ij} = \tilde{\omega}_j -\tilde{\omega}_i.\eqno(6)
$$
Note that this implies that the map $H^1(B_2) \to H^1(Z_1)$ has a  non-trivial 
kernel. Suppose that the left-hand-side of (6) represents an element in the
image of $H^1(B_1)=H^1(dO_X) \to H^1(B_2)$. Then we find a relation $
dh_{ij}= \hat{\omega}_j - \hat{\omega}_i$ with $\hat{\omega}_i$ closed. Then
since $C$ annihilates $dh_{ij}$ the $C\hat{\omega}_i$ define a  global
$1$-form and this must be zero. Hence we can write $\hat{\omega}_i=d\phi_i$
and this shows that $dh_{ij}$ represents the zero-class in $H^1(dO_X)$,
contrary to the assumption. Hence we find a  non-trivial element in
$H^1(B_2)$ which does not lie in the image of the natural inclusion
$H^1(B_1) \to H^1(B_2)$. Carrying out  this argument for all $n$ proves the
claim. 
$\square$
\par
\proclaim (\chno.2) Corollary. Assume that  $h < \infty$.  Then for all $n\geq
1$  the natural map $\beta_n:H^1(B_n) \to H^1(X, \Omega^{1}_{X})$ is injective
and the image has dimension $\min\{n,h-1\}$. \par
\noindent
{\sl Proof.}   Since by (6.2) $\dim H^1(B_n)$ stabilizes for $h\neq \infty$,
non-injectivity would contradict the preceding proposition. $\square$
\par
Note that the natural map $H^1(B_n) \tto H^1(\Omega^{1}_{X})$ 
is not necessarily injective
for $h=\infty$ because $\dim H^1(B_n) > 20$ for $n>20$.
In the case of $h \neq \infty$, we often identify 
$H^{1}(B_{n})$ with the image of the natural inclusion 
$H^{1}(B_{n}) \to H^{1}(\Omega_{X}^{1})$ in Corollary 9.2.
\bigskip
Let $Z_n \to \Omega^{1}_{X}$ be the natural inclusion. We have an induced map
 $$
z_n: H^1(Z_n) \tto H^1(\Omega^{1}_{X}).
$$ 
We would like to characterize both the image and the kernel of this map.
We often write ${\rm Im}(H^1(Z_n))$ for the image of $z_n$.
\par
\proclaim (\chno.3) Lemma. i) We have ${\rm Im}(H^1(B_n))\subseteq {\rm
Im}(H^1(Z_n))^{\bot}$, in particular, ${\rm Im}(H^1(B_n))\subseteq {\rm
Im}(H^1(B_n))^{\bot}$.  ii) Assume that $h< \infty$.  If $C: H^1(B_{n+1})
\to H^1(B_n)$ is surjective then we have the equality
${\rm Im}(H^1(Z_n)) = {\rm Im}(H^1(X,B_n))^{\bot}$. \par
\noindent
{\sl Proof.} We first show that ${\rm Im}(H^1(B_n))$ and ${\rm Im}(H^1(Z_n))$
are orthogonal.  Let $\alpha\in {\rm Im}(H^1(B_1))$ and $\beta
\in {\rm Im}(H^1(Z_1))$. Then we find an element $\alpha \wedge \beta \in
H^2(\Omega^{2}_{X})\cong k$ representing the cup product $\langle \alpha, \beta
\rangle$. If we apply Cartier  to $\alpha \wedge \beta$ 
a suitable number of times then it is zero. Now use
the exact sequence
$$
0 \to d\Omega^{1}_{X} \tto \Omega^2_{X, \rm closed} {\buildrel C \over \tto}
\Omega^{2}_{X} 
\to 0, \eqno(7)
$$
and the fact that $\Omega_{X, \rm closed}^2 = \Omega^2_X$. Then,
we have from the the long exact sequence the exact sequence
$$
 H^2(d\Omega^1_X) \tto H^2(\Omega^2_X) {\buildrel C \over \tto} H^2(\Omega^2_X)
\to 0. 
$$
The fact that $\dim H^2(\Omega^2_X)= 1$ implies that 
$C: H^2(\Omega^2_X)\to H^2(\Omega^2_X)$ is an isomorphism as a $p$-linear
mapping.  Therefore, 
for $x \in H^2(\Omega^2_X)$ we have $x=0$ if and only if $C^n(x)=0$ 
for some $n$. Hence, we conclude 
$\alpha \wedge \beta = 0$.
 We now prove equality by induction. For $n=1$ we have
${\rm Im}(H^1(B_1))^{\bot}= {\rm Im}(H^1(Z_1))$ because ${\rm
Im}(H^1(B_1))$ is the kernel $G_1/F^2$ of the Cartier operator and ${\rm
Im}(H^1(Z_1))$ is $F^1\cap G_2= F^1\cap G_1^{\bot}$. Suppose that we have
proved that ${\rm Im}(H^1(B_i))^{\bot}= {\rm Im}(H^1(Z_i))$ for $i\leq n$.
If $\beta \in {\rm Im}(H^1(Z_1))$ is orthogonal to all  $\alpha \in {\rm
Im}(H^1(B_{n+1}))$ then we have $\langle C\alpha, C\beta\rangle =0$ and
since $C: H^{1}(B_{n+1}) \to H^{1}(B_{n}) $ is surjective this implies 
that $C\beta \in {\rm Im}(H^1(Z_{n}))$, i.e.\
$\beta \in {\rm Im }(H^1(Z_{n+1}))$. $\square$
\smallskip 
\noindent
\proclaim (\chno.4) Lemma. The Cartier operator $C: H^1(B_n) \to H^1(B_{n-1})$
is surjective for $n\leq h-1$. Moreover, for $n\leq h-1 < \infty$ we have 
$\dim {\rm Im}(H^1(Z_n)) = 20-n$
\par
\noindent
{\sl Proof.} Note that we know that $h^1(B_n)=n$ for $n\leq h-1$
and thus the exact sequence $0\to B_1\to B_n \to B_{n-1} \to 0$
implies that $C: H^1(B_{n})\to H^1(B_{n-1})$ is surjective for
$n\leq h-1$. The rest follows from (9.3).~$\square$
\par
\smallskip
\proclaim (\chno.5) Corollary. If $h\neq \infty$ we have the following orthogonal
filtration in $H^1(\Omega^1)$:
$$
\eqalign { 0\subset H^1(B_1) \subset H^1(B_2) \subset \ldots & \subset
H^1(B_{h-1}) \subset 
\cr
\subset {\rm Im}(H^1(Z_{h-1})) &\subset {\rm Im}(H^1(Z_{h-2} )) \subset \ldots
\subset {\rm Im}(H^1(Z_1)) \subset H^1(\Omega^1_X). \cr} \eqno(8)
$$
\par
The exact sequence (3) gives for $i=0$ rise to the exact sequence
$$
0 \to H^0(d\Omega^1_X) \tto H^1(Z_1) \tto H^1(\Omega^1_X)
{\buildrel d \over\tto} H^{1}(d\Omega^1_X) \tto H^2(Z_1) \to
$$
The natural map $H^1(Z_1) \to H^1(\Omega^1)$ is the composition of
$H^1(Z_1)\to H^2_{\rm dR}$ and the projection $H^2_{\rm dR}\to F^1/F^2$, i.e.\ by
(7.1) it is the map $F^1\cap G_2 \to F^1/F_2$. This is an isomorphism for $h=1$
and it has a $1$-dimensional kernel otherwise.  It follows that 
$$
\dim H^0(d\Omega^1)=\dim H^1(d\Omega^1)=\cases{0 &if $h=1$\cr
1& if $h\neq 1$.\cr}
$$
From the exact sequence 
$$
0\to H^0(d\Omega^1_X) \tto H^1(Z_{n+1}){\buildrel \psi_{n+1} \over\ttto}  
H^1(Z_n)
\tto H^1(d\Omega^1_X) \to 
$$
with $\psi_{n+1}$ the map induced by inclusion we deduce that for $h\neq 1$
$$
\psi_{n+1} \hbox{ \rm is  surjective } \iff \dim H^1(Z_{n+1}) > \dim H^1(Z_n).
\eqno (9)
$$
\proclaim (\chno.6) Lemma. For $h\neq \infty$ we have $\dim
H^1(Z_n)=20$. \par
\noindent
{\sl Proof.} If $h=1$ we have $h^0(d\Omega^1_X)=0$ and $h^1(d\Omega^1_X)=0$  
hence all $\psi_n$ are isomorphisms. Since we know $h^1(Z_0)=20$ the result
follows for $h=1$. If $h\neq 1$ then $H^0(d\Omega^1_X)\cong k$.  For $n \leq
h-1$ we have ${\rm Im}(H^1(Z_n)) = 20-n$ by Lemma 9.3
and $\dim H^{1}(B_{n}) = \min \{n, h-1\}$. Suppose there exists an  $n ~(n \leq
h - 1)$  such that $\psi_{n}$ is surjective.
Take the smallest such $n$. Then, the image of $z_{n}$
coincides with the image of $z_{n-1}$, which
contradicts $\dim {\rm Im}(H^1(Z_n)) \neq  {\rm Im}(H^1(Z_{n-1}))$. Hence, 
$h^1(Z_n) =20$ for $n\leq h-1$. Consider for $n=h$ the commutative diagram of
exact sequences
$$
\def\normalbaselines{\baselineskip20pt
\lineskip3pt \lineskiplimit3pt }
\matrix
{&&0&&0&&&&  \cr
  &  &\mapdown{}&  & \mapdown{} && \cr&&B_1& = &B_1&&&&  \cr
  &  &\mapdown{}&  & \mapdown{} && \cr
 0&\tto&B_n&\tto&Z_n&{\buildrel C^n \over \tto}&\Omega^1_X &\tto&0.\cr
  &  &\mapdown{C}&  &\mapdown{C} &&\mapdown{=} \cr
 0&\tto&B_{n-1}&\tto&Z_{n-1}&{\buildrel C^{n-1} \over \tto}&\Omega^1_X &\tto&0.\cr
  &  &\mapdown{}&  &\mapdown{} && \cr
 &&0&&0&&&&.\cr
}
$$
The diagram shows that $C^{h}: H^1(Z_h) \tto H^1(\Omega^1_X)$ 
factors through the image of $C^{h-1}$. This implies that
$\dim H^1(Z_h) - (h-1)\leq 20 - (h-1)$. Since $h^1(Z_n)\geq 20$ for all
$n\geq 1$ we get $h^1(Z_h)=20$. We can repeat this argument for
$H^1(Z_m)$ with $m\geq h$.
$\square$ 
\bigskip
\noindent
\def\chno{10}\centerline{\bf \chno. Chern classes of line bundles and  closed
forms}
\smallskip
\noindent
We start with a well-known result due to Ogus [O, Cor.\ 1.5]. 
We give here the proof by Shafarevich [Sh] for the reader's
convenience.
\smallskip
\proclaim (\chno.1) Proposition. 
The map $c_1: NS(X)/pNS(X) \tto
H^2_{\rm dR}$ is injective and factors through $F^1H^2_{\rm dR}$. \par
\noindent
{\sl Proof.} (Shafarevich) We take an affine open covering $\{U_i\}$ of $X$.
A class in $H^2_{\rm dR}$ is represented by a tripel $(a,b,c)\in
 C^2(O_X) \oplus C^1(\Omega^1_X) \oplus C^0(\Omega^2_X)$. The boundaries are of
the form $(\delta h_{ij}, dh_{ij}+\omega_j -\omega_i, d\omega_i)$ with
$(h_{ij},\omega_i) \in C^1(O_X) \oplus C^0(\Omega^1_X)$.  So if a Chern class
$c_1(L)$, represented by $(0, d\log f_{ij},0)$, is zero in $H^2_{\rm dR}$
then  there exists $(h_{ij},\omega_i) \in C^1(O_X)\oplus C^0(\Omega^1_X)$
with $d\omega_i=0$ and
$\delta h_{ij}=0$ and we have $ d\log f_{ij}= \omega_j - \omega_i + dh_{ij}$. 
By the relation $\delta(h_{ij})=0$ the $h_{ij}$ defines a class in
$H^1(X,O_X)=0$, so we have $h_{ij}=\eta_j-\eta_i$ with $\eta_i$ regular  and
we can replace
$\omega_i$ by $\omega_i+d\eta_i$ and obtain a relation
$$
d\log f_{ij}= \omega_j - \omega_i \qquad {\rm with} \quad \omega_i \quad {\rm
closed}. \eqno(10)
$$
Applying the Cartier operator we find
$$
d\log f_{ij}= C\omega_j - C \omega_i. \eqno(11)
$$
subtracting (1) from (2) we find $C\omega_i-\omega_i=C\omega_j-\omega_j$.
This defines a global $1$-form which must be zero. Hence we see 
$C\omega_i=\omega_i$ and it follows that $\omega_i= d\log \phi_i$ (after
shrinking the $U_i$ if necessary). We find
$$
d\log f_{ij} = d\log \phi_i \phi_j^{-1},
$$
hence
$$
f_{ij}=  \phi_i \phi_j^{-1} \psi_{ij}^p.
$$
for some $\psi_{ij} \in O(U_{i} \cap U_{j})$. Thus modulo a $p$-th power 
$L$ is trivial. The proof also shows that the image lands in $F^1H^2_{\rm dR}$.
$\square$ 
\par
\proclaim (\chno.2) Proposition. If $h< \infty$ then we  have $\langle
c_1(NS(X))\rangle \cap {\rm Im}(H^1(B_n)) = \{ 0 \}$ for all~$n$.
Moreover, $c_1(NS(X))$ is orthogonal with ${\rm Im}(H^1(B_n))$ for all $n$.\par
\noindent
{\sl Proof.}  First we show that $c_1(NS(X)) \cap {\rm Im}(H^1(B_n))= (0)$
for all $n>0$. If it is not, then take a minimal $n$ such that ${\rm
Im}(H^1(B_n))$ contains a Chern class $ 0 \neq [d\log f_{ij}]$. 
We can write a (non-trivial) relation as
$$
d\log f_{ij}=\beta_{ij} + \omega_j -  \omega_i,\eqno(12)
$$
where the $\beta_{ij}$  are forms in
$B_n$,  but not in $B_{n-1}$. Apply the inverse Cartier operator as in (9.1) to
get a relation
$$
d\log f_{ij} = \tilde{\beta}_{ij} + \tilde{\omega}_j - 
\tilde{\omega}_i\eqno(13)
$$   
where the $\tilde{\omega}_i$ are closed forms with
$C(\tilde{\omega}_i)=\omega_i$
 and the $\tilde{\beta}_{ij}$  are forms in $B_{n+1}$ with
$C(\tilde{\beta}_{ij})= \beta_{ij}$. Subtracting (12) from (13) shows that
$\tilde{\beta}_{ij}-{\beta}_{ij}$ is a  boundary.
Since $\beta_{ij} $ defines a non-zero element of $H^1(B_n)$ 
which is not in the image of $H^1(B_{n-1})$ the cocycle 
$\tilde{\beta}_{ij}$ gives an element of $H^1(B_{n+1})$ not 
in the image of $H^1(B_n)$. Hence the left hand side is not zero in
$H^1(B_{n+1})$ and this shows that $H^1(B_{n+1})\to H^1(\Omega^1)$ is not
injective.
\par
Suppose now that $\langle c_1(NS(X))\rangle \cap {\rm Im}(H^1(B_n))\neq
0$. Considering all $n$ which satisfy this condition, we then have a relation 
with $m\geq 2$ minimal
$$
d\log f_{ij}^{(1)} + \sum_{\nu=2}^m a_{\nu} d\log f_{ij}^{(\nu)} =
\beta_{ij} + \omega_j -\omega_i.
$$
We may assume that $m\geq 2$ and that $a_{\nu}\not\in \F_p$ for all $\nu
\geq 2$. Then by applying
$C^{-1}$ as before we find 
$$
d\log f_{ij}^{(1)} + \sum_{\nu=2}^m a_{\nu}^{p} d\log f_{ij}^{(\nu)} =
\tilde{\beta}_{ij} +
\tilde{\omega}_j -\tilde{\omega}_i,
$$
where the $\tilde{\omega}_i$ are closed and $C\tilde{\omega}_i=\omega_i$.
Subtracting the two relations we  find a shorter relation ($m$
smaller but with $n$ maybe larger). This contradiction shows that $\langle
c_1(NS(X))\rangle
\cap {\rm Im}(H^1(B_n))= 0$. 
\par
The orthogonality of $\langle c_1(NS(X))\rangle$ and ${\rm Im}(H^1(B_n))$
follows from the fact that $\langle c_1(NS(X))\rangle \subset {\rm
Im}(H^1(Z_n))$ and Lemma (9.3).  
$\square$
\par
\proclaim (\chno.3) Proposition. Suppose that  $h< \infty$. Then the Chern class
map $$c_1 \otimes k: NS(X)/pNS(X)\otimes k \to H^1(X, \Omega^1_X)$$ is injective.
\par
\noindent
{\sl Proof.} Suppose we have a relation
$\sum_{\nu =1}^r a_\nu c_1(L_\nu ) =
\omega_j - \omega_i$ for line bundles $L_\nu $ and $a_\nu \in k$. We may assume
that
$a_1=1$ and that the relation is the shortest possible ($r$ minimal).
Furthermore, we can assume that $a_\nu /a_\mu \not\in \F_p$ for  $\nu \neq \mu$;
otherwise we can easily find a shorter one. Now apply the inverse Cartier
operator $C^{-1}$ to the relation as we did before. We find a new relation
$$
dg_{ij}+ c_1(L_1) + \sum_{\nu =2}^r a_\nu^p c_1(L_\nu )- \tilde{\omega}_i +
\tilde{\omega}_j =0, 
$$
where the $g_{ij}$ are regular on $U_i \cap U_j$. If the cocycle
$dg_{ij}$ defines a zero class in $H^1(X, \Omega^1_X)$,  we can write
$dg_{ij}= \eta_j -
\eta_i$, and we can replace the relation by a shorter one by subtracting the
two  relations contradicting the minimality of $r$. Hence $\{ dg_{ij}\}$ defines
a non-zero class in $H^{1}(X, \Omega^{1}_{X})$ and we find  a non-zero element 
in ${\rm Im}(H^1(B_1)) \cap \langle
c_1(NS(X))
\rangle $. $\square$
\smallskip
\par
As a corollary of (6.2), (10.2) and (10.3) we now find the well-known result of
Artin and Mazur on the rank $\rho$ of the N\'eron-Severi group:
\par
\proclaim (\chno.4) Corollary. For $h\neq \infty$ we have $\rho \leq 22-2h$. 
\par
\noindent
{\bf (\chno.5) Remark.} A line bundle $L$ defined by transition functions
$f_{ij}$ defines a cocycle
$d\log f_{ij}$ with values in $Z_n\Omega^1_X$ for all $n\geq 0$. We thus can
view the class $c_1(L)$ as a class in $H^1(Z_n)$ for all $n\geq 0$ as well as
in $H^2_{\rm dR}$. If $h<\infty$ the maps 
$$
c_1\otimes k: NS(X)/pNS(X) \otimes k \tto H^1(Z_n)
$$
are injective for all $n\geq 0$. 
\bigskip
\noindent
\def\chno{11}\centerline{\bf \chno. The Supersingular Case}
\smallskip
\noindent
The map $c_1: NS/pNS \to H_{\rm dR}^2$ is injective and factors through
$H^1(Z_j)$ for all $j\geq 1$. However, the map $c_1\otimes k: NS \otimes k
\to H^2_{\rm dR}$ is not necessarily injective. For $X$ supersingular in
Shioda's sense, i.e.\ $\rho=B_2=22$, it cannot be injective since $\dim_k
H^1(Z_1)=20$ or $21$, the latter if $X$ is superspecial. 
\smallskip
We define for $j=0,1\ldots$
$$
U_j := \ker\{ c_1\otimes k: NS\otimes k \to H^1(Z_j)\}
$$
and we set 
$$
\dim U_1=\sigma_0.
$$  
Using  the natural maps $H^1(Z_i) \to
H^1(Z_{i-1})$ we have $ U_{j+1} \subseteq U_{j}$ for  $j=0,1,2,\ldots $.
We define two bijective operators  on $NS\otimes k$
$$
\varphi=1\otimes F \qquad {\rm  and } \qquad 
\gamma= 1\otimes F^{-1},
$$ 
with the Frobenius action $F: a\mapsto a^p$ on the second factor $k$. 
\smallskip
\noindent
{\bf (11.0) Remark.} If we assume that $\rho=B_2=22$ (i.e.\ the truth of the Artin
conjecture that $h=\infty$ implies $\rho=22$) then one can show that the invariant
$\sigma_0$ just introduced equals the Artin invariant $\sigma_0$, i.e.\ the
intersection form on the lattice
$NS(X)$ has discriminant 
$$
{\rm disc}(NS(X))=-p^{2\sigma_0}.
$$
\par
\proclaim (\chno.1) Lemma. We have $\gamma(U_{j+1}) \subseteq U_j$;
equivalently, we have $ U_{j+1}\subseteq \varphi(U_j) $. Moreover, we have 
$U_{j+1} \subseteq U_j \cap \varphi(U_{j})$.\par
\noindent
{\sl Proof.}  This follows from the  commutativity of the
diagram
$$
\def\normalbaselines{\baselineskip20pt
\lineskip3pt \lineskiplimit3pt }
\matrix
{NS\otimes k &{\buildrel \gamma \over \tto} &NS \otimes k &&&&&&   \cr
 \mapdown{c_1\otimes k} &              & \mapdown{c_1\otimes k}&  &
&&\cr
 H^1(Z_{j+1}) &{\buildrel C \over \tto}&H^1(Z_{j})&&&&&&\cr}
$$
with $C$ the Cartier operator. The second result follows from this and the
inclusion $U_{j+1} \subset U_j$.
$\square$
\par
Now choose an element $u_{\min}= u_{\min}^{(j)} \neq 0$ of minimal length in
$U_j$ under  the assumption that $U_j$ is non-zero, i.e.\ write 
$u_{\min} = \sum_{i=1}^m a_i [L_i]$ and require $m\geq 2$ to be minimal.  We
also may assume --and we shall-- that $a_1=1$.
\bigskip
\proclaim (\chno.2) Lemma. For $j\geq 1$ we have $u_{\min} \not\in
\varphi(U_{j})$.   Similarly we have
$u_{\min}\not\in~\gamma(U_j)$.  If $X$ is not superspecial the
conclusion holds also for $j=0$. \par
\noindent
{\sl Proof.}  If $u_{\min}\in \varphi(U_j)$ is such a minimal  element with
$a_1=1$ then $\gamma(u_{\min}) -u_{\min}$ would be a shorter element or zero. If
it is zero, then $u_{\min} \in NS\otimes {\bf F}_p
\cap U_1 = \{ 0 \}$. For $j=0$ the argument is similar. Note that 
$ NS\otimes {\bf F}_p \cap U_{0} \neq \{0 \}$ if and only if $X$ is 
superspecial, cf.\ [O 1], Cor. 1.4.
\par
\proclaim (\chno.3) Lemma. The map $c_1\otimes k : \varphi(U_j) \to
H^1(Z_{j+1})$ factors via $H^1(B_1) \to H^1(Z_{j+1})$ and the induced  map
$\varphi(U_j)  \to H^1(B_1)$ is surjective if $U_j \neq \{ 0 \}$. \par
\noindent
{\sl Proof.} If $u \in U_j$ there exist closed forms $\zeta_x \in Z_j(V_x)$
for some open covering $V_x$ such that $(c_1\otimes k)(u)$ is a coboundary: 
$\zeta_{\beta} - \zeta_{\alpha}$. Now use the local surjectivity of $C$ to write
$$
(c_1\otimes k)(\varphi(u))= \tilde{\zeta}_{\beta} -\tilde{\zeta}_{\alpha} + 
\phi_{\alpha \beta}
$$
with $\tilde{\zeta}_{x} \in Z_{j+1}$, $C\tilde{\zeta}_{x} = \zeta_{x}$, 
$\phi_{\alpha \beta} \in B_1$ on a suitable open covering. Then this
$\phi_{\alpha \beta}$ defines a cocycle, thus an element in $H^1(B_1) \subset
H^1(Z_{j+1})$.
\par
To prove the surjectivity, choose a non-zero element $u_{\min} \in U_j$. 
Suppose that $\phi_{\alpha \beta} = \eta_{\beta}-\eta_{\alpha}$ with $\eta \in
B_1$. Then $\varphi (u_{\min}) \in U_j$, hence $U_{\min} \in \gamma(U_j)$ which
contradicts Lemma (\chno.2).
\par
\proclaim (\chno.4)  Corollary. We have $U_{j+1} = U_j \cap \varphi(U_j)$ 
and $\dim (U_{j+1})=\max\{ \dim(U_j)-1,0\}$.\par
\noindent
{\sl Proof.}  The kernel of $c_1\otimes k : \varphi(U_j) \to H^1(Z_{j+1})$ equals
$U_{j+1}$ by (11.1) and has codimension $1$ by (11.3). Since  $U_j \neq
\varphi(U_j)$, and since their intersection contains $U_{j+1}$ we must have 
$U_{j+1} = U_j \cap \varphi(U_j)$. The statement about dimensions follows.
\par
If we assume that $\sigma_0\geq 1$ then we have a strictly increasing sequence
$$
\{ 0 \} = U_{\sigma_0+1}  \subset U_{\sigma_0} \subset \ldots \subset U_2
\subset U_1 \eqno(14)
$$
and this implies:
\par
\proclaim (\chno.5) Proposition. The map $c_1\otimes k$ factors through an 
{\sl injection } $$ NS(X)/pNS(X)\otimes k \to H^1(Z_{\sigma_0+1}).$$
\par
We can generalize the result of Corollary (11.4).
\par
\proclaim (\chno.6) Lemma. We
have
$\varphi^k(U_j)\cap U_j= U_{j+k}$. In particular $\varphi^{\sigma_0}(U_1)\cap
U_1 =\{ 0 \}$.
\par
\noindent
{\sl Proof.} We prove this by induction on $k$, the case $k=1$ 
was proved in (\chno.4). Suppose it holds for $k$. Then 
$$
\varphi^{k+1}(U_j) \cap U_j \subset \varphi[ \varphi^k(U_{j-1}) \cap
U_{j-1}] \subset \varphi(U_{j+k-1})
$$
On the other hand we have
$$
\varphi(U_{j+k-1})\cap \varphi^{k+1}(U_j) \subset
\varphi(U_j \cap \varphi^k(U_j))= \varphi(U_{j+k}).
$$
But by an easy induction one has
$$
\varphi(U_{j+k})\cap U_j \subset U_{j+k+1}.
$$
In view of $\dim (\varphi^{k+1}(U_j) \cap U_j)\geq \dim(U_j)-(k+1)$ the
result follows.
\proclaim (\chno.7) Lemma. Suppose that $U_1\neq \{ 0\}$ and let $u_{\min} \in
U_1$. Then $\gamma(u_{\min})\in U_0\backslash U_1$. In particular, $(c_1\otimes
k)(\gamma(u_{\min}))\in H^0(\Omega^2)\subset H^1(Z_1)\subset H^2_{\rm dR}$. \par
\noindent
{\sl Proof.} Since $\gamma(u_{\min})$ does not lie in  $U_1$, but lies in $U_0$
we see that $(c_{1}\otimes k)(\gamma(u_{\min}))$ must lie in the kernel of
$H^1(Z_1)
\to H^1(\Omega^1)$,  which  is
$H^0(\Omega^2_X)$.\par
\proclaim (\chno.8) Lemma. The Chern class map $c_1\otimes k :
\varphi^m(U_j) \tto H^1(Z_{j+m})$ factors through $H^1(B_m)$. For any $t\geq 1$ 
the natural image of $H^1(B_t)$ in $H^1(Z_{\sigma_0+1})$ is contained in
the image of $NS(X)/pNS(X) \otimes k$ under $c_1\otimes k$. 
\par
\noindent
{\sl Proof.} As in the proof of (11.3) we can write
$(c_1\otimes k)(u) = \zeta_{\beta} - \zeta_{\alpha}$ with $\zeta_{x} \in
Z_{j}(V_{x})$. Now use the local surjectivity of $C$ to write
$$
\varphi^m(u)= \tilde{\zeta}_{\beta} -\tilde{\zeta}_{\alpha} + \phi_{\alpha
\beta}
$$
with $\tilde{\zeta} \in Z_{j+m}$, $C^m\tilde{\zeta}_{x} = \zeta_{x}$,
$\phi_{\alpha
\beta} \in B_m$. Then this $\phi_{\alpha \beta}$ defines a cocycle, thus an
element in $H^1(B_m) \subset H^1(Z_{j+m})$. This proves the first statement.
\par
We prove the second statement by induction. Note that by (11.3)  the image of
$H^1(B_1)$ in $H^1(Z_{\sigma_0+1})$ is contained in the image of $NS(X)/pNS(X)
\otimes k$ under $c_1\otimes k$. Let $\alpha$ be an element of the image of 
$H^1(B_t)$ and $\beta= C\alpha$ in the image of $H^1(B_{t-1})$. Then $\beta =
(c_1\otimes k)(v)$ for some $v \in NS\otimes k$. But then $\alpha -(c_1\otimes
k)(\varphi(v))$ is an element of $H^1(B_1)$. By induction this is in the image of
$(c_1\otimes k)(NS \otimes k)$. Hence $\alpha$ lies in the image of
$(c_1\otimes k)(NS\otimes k)$.
\par
\proclaim (\chno.9) Proposition. Let $\sigma_0 \geq 1$. The dimension of  the
image of
$H^1(B_{\sigma_0})$ in
$H^1(Z_1)$ equals $\sigma_0$. The image in $H^1(\Omega^1_X)$ is
$\sigma_0-1$-dimensional. \par
\noindent
{\sl Proof.} The first statement follows directly from (\chno.6) and (\chno.8).
Arguing similarly for $U_0$ we find that $c_1\otimes k:\varphi^{\sigma_0}(U_0) 
\to H^1(\Omega^1_X)$ factors through the natural map $H^1(B_{\sigma_0})\to
H^1(\Omega^1_X)$. The intersection $\varphi^{\sigma_0}(U_0)\cap U_0$ has
dimension $1$.
\par
\proclaim (\chno.10) Theorem. For a K3 surface $X$ with $B_{2} = \rho$ and Artin
invariant $\sigma_0$,  we have $\dim ({\rm Im}~H^{1}(Z_{\sigma_0}))
= 21 - \sigma_0$ for the image in $H^1(\Omega^1_X)$ and it is generated by Chern
classes. \par
\noindent
{\sl Proof.} Since we have
$$
\langle c_{1}(NS(X)/pNS(X)) \rangle \subset {\rm Im}~ H^{1}(Z_{\sigma_0})
\subset ({\rm Im}~H^{1}(B_{\sigma_0}))^{\perp}
\subset H^{1}(\Omega^{1}_{X})
$$
and $\dim \langle c_{1}(NS(X)/pNS(X)) \rangle =
\dim ({\rm Im}~ H^{1}(B_{\sigma_0}))^{\bot} = 20 - (\sigma_0 - 1)$ by (\chno.9),
we have 
$$
\langle c_{1}(NS(X)/pNS(X)) \rangle =
{\rm Im}~H^{1}(Z_{\sigma_0}) = ({\rm Im}~H^{1}(B_{\sigma_0}))^{\bot}$$
and so $\dim{\rm Im}~H^{1}(Z_{\sigma_0}) = 21- \sigma_0$. $\square$
\smallskip
Since the codimension of ${\rm Im}~ H^{1}(Z_{i+1})$ in ${\rm Im} H^{1}(Z_{i})$
is at most one, we conclude that
$$
{\rm Im}~H^{1}(Z_{\sigma_0}) = {\rm Im}~ H^{1}(Z_{\sigma_0 -1}) \subset
{\rm Im}~ H^{1}(Z_{\sigma_0 -2}) \subset \ldots \subset {\rm Im}~H^{1}(Z_{1})
\subset  H^{1}(\Omega^{1}_{X})
$$
and ${\rm Im}~H^{1}(Z_{n}) = {\rm Im}~H^{1}(Z_{\sigma_0})$ for $n \geq \sigma_0$.
Here, the inclusions are strict inclusions. Moreover, we see that
the injection
$c_{1} \otimes k : NS(X)/pNS(X) \otimes k \to H^{1}(Z_{\sigma_0 +1})$
is an isomorphism:
$$
c_{1} \otimes k : NS(X)/pNS(X) \otimes k \cong  H^{1}(Z_{\sigma_0 +1}).
$$
\par
We now need the following lemma.
\par
\proclaim (\chno.11) Lemma. Let $X$ be a K3 surface $X$ with $B_{2} = \rho$  and
Artin invariant $\sigma_0$. For every $n\geq 0$ the natural map
$H^{1}(Z_{\sigma_0 + n+ 1}) \to H^{1}(Z_{\sigma_0 +n})$ is surjective. \par
\noindent
{\sl Proof.}  By Theorem (11.10) the dimension of the image of
$H^1(Z_{\sigma_0})$ in $H^1(\Omega^1)$ is
$21-\sigma_0$. By (14) it follows that the image of $H^1(Z_{\sigma_0-1})$ in
$H^1(\Omega^1)$ has dimension at least $22-\sigma_0-1$. Since the map
$H^1(Z_{\sigma_0+1}) \to H^1(\Omega^1)$ factors through $H^1(Z_{\sigma_0})$
the map $H^1(Z_{\sigma_0+1}) \to H^1(Z_{\sigma_0})$ must be surjective.
\par
We now prove that if the natural mapping $H^{1}(Z_{n+1}) \to
H^{1}(Z_{n})$ is surjective, then so is $H^{1}(Z_{m+1})  \to H^{1}(Z_{m})$ for
any $m \geq n$. Suppose that the natural homomorphism 
$H^{1}(A, Z_{n+1}) \rightarrow H^{1}(A, Z_{n})$ is surjective.
By the diagram of exact sequences
$$
\def\normalbaselines{\baselineskip20pt
\lineskip3pt \lineskiplimit3pt }
\matrix
{0&\to &B_1& \longrightarrow &Z_{n+2}&{\buildrel C \over \longrightarrow}
&Z_{n+1}& \to 0 \cr
  && \mapdown{=} &&\mapdown{\iota_{n+2}}&  & \mapdown{\iota_{n+1}}  \cr
0&\to &B_1& \longrightarrow &Z_{n+1}&{\buildrel C \over \longrightarrow}&
Z_{n}& \to 0 \cr
}
$$
we have a diagram of  exact sequences
$$
\def\normalbaselines{\baselineskip20pt
\lineskip3pt \lineskiplimit3pt }
\matrix
{\rightarrow & H^{1}(X, B_{1}) & \rightarrow &
H^{1}(X, Z_{n + 2}) & {\buildrel C \over \longrightarrow} & 
H^{1}(X, Z_{n+1}) & \rightarrow &
H^{2}(X, B_{1}) \cr
  & \mapdown{=} &  & \mapdown{\iota_{n+2}} &    &\mapdown{\iota_{n+1}} &  
  & \mapdown{=}   \cr
      \rightarrow & H^{1}(X, B_{1}) & \rightarrow &
H^{1}(X, Z_{n +1}) & {\buildrel C\over \longrightarrow} & H^{1}(X, Z_{n}) & 
\rightarrow & H^{2}(X, B_{1}).\cr
}
$$
From this diagram  we see that the natural homomorphism 
$H^{1}(X, Z_{n+2}) \rightarrow H^{1}(X, Z_{n+1})$
is also surjective. So this lemma now follows by induction. $\square$
\par
\proclaim (\chno.12) Corollary. Let $X$ be a K3 surface $X$ with $B_{2} =  \rho$
and Artin invariant $\sigma_0$. For $n\geq \sigma_0$  we have ${\rm
Im}(H^1(B_n))={\rm Im}H^1(Z_n)^{\bot}$ and  $\dim{\rm Im}~H^{1}(B_{n}) = \sigma_0
- 1$.
\par
\noindent
{\sl Proof.} By the proof of (11.10), we have
${\rm Im}H^{1}(Z_{\sigma_{0}})^{\perp} = {\rm Im}H^{1}(B_{\sigma_{0}})$.
Therefore, for $n \geq \sigma_{0}$, we have
$$
{\rm Im}H^{1}(Z_{n})^{\perp} = {\rm Im}H^{1}(Z_{\sigma_{0}})^{\perp} =
{\rm Im}H^{1}(B_{\sigma_{0}}) \subset {\rm Im}H^{1}(B_{n}).
$$
On the other hand, by the proof of (9.3), we have
${\rm Im}H^{1}(Z_{n})^{\perp} \supset {\rm Im} H^{1}(B_{n})$.
Hence, we get the desired results.
 $\square$
\par
Since $c_{1} \otimes k : NS(X)/pNS(X) \otimes k \tto  H^{1}(Z_{i})$
is injective for $i \geq \sigma_0 + 1$, we have the following proposition.
\par
\proclaim (\chno.13) Proposition. For a K3 surface $X$ with $B_2 = \rho$ the
following four conditions are equivalent.
\item{(i)} The natural map $H^{1}(Z_{i}) \to H^{1}(Z_{i-1})$ is surjective.
\item{(ii)} The Cartier operator $C : H^{1}(Z_{i}) \to H^{1}(Z_{i-1})$ is
surjective.
\item{(iii)} $\dim H^{1}(Z_{10}) \geq 31 - i$.
\item{(iv)} $\sigma_0 \leq i$. \par
\bigskip
\bigskip
\noindent
\def\chno{12}\centerline{\bf \chno. The Kodaira-Spencer Map}
\smallskip
\noindent
Let $X_{0}$ be a K3 surface, and let $\pi : X \tto S$  be the versal formal
k-deformation of $X_{0}$.  Then, as is well-known (cf.\ [D]), we have $S = {\rm
Spf} k[[t_{1}, \ldots, t_{20}]]$  with variables $t_{1}, \ldots, t_{20}$. 
We denote by $\nabla$ the Gauss-Manin connection of
$H^{2}_{\rm dR}(X/S)$:
$$
\nabla : H^{2}_{\rm dR}(X/S) \tto 
\Omega_{S/k}^{1}\otimes H^{2}_{\rm dR}(X/S).
$$
We take a basis $\omega$ of
$H^{0}(X, \Omega^{2}_{X/S})$. Then, $\nabla$ composed  with cup
product with
$\omega$ gives an isomorphism:
$$
\rho_{\omega} : H^{1}(X, \Omega^{1}_{X/S}) {\buildrel \sim \over \tto }
 \Omega^{1}_{S/k}.
$$
We denote by $m$ the maximal ideal of the closed 
point of $S$. By evaluating $\rho_{\omega}$ at zero
we have an isomorphism:
$$
\rho_{\omega,0} : H^{1}(X_{0}, \Omega^{1}_{X_{0}/k}) 
{\buildrel \sim \over \tto } m/m^{2}.
$$
{\bf (\chno.1) Remark.}  Ogus gave an  explicit expression of the
isomorphism $\rho_{\omega}$ as follows. For an element
$\alpha \in  H^{1}(X, \Omega^{1}_{X/S})$ we choose 
a lifting $\alpha^{\prime} \in F^{1}H^{2}_{\rm dR}(X/S)$ of
$\alpha$. Since 
$\langle \alpha^{\prime}, \omega\rangle = 0$, we have
$$
     \rho_{\omega}(\alpha) 
= \langle \nabla \alpha^{\prime},  \omega \rangle
= - \langle \alpha^{\prime}, \nabla \omega \rangle.
$$
For details, see the paper by Deligne/Illusie [D], cf. also Ogus $[O]$.
\bigskip
\def\chno{13}\centerline{\bf \chno.  Horizontality }
\smallskip 
\noindent
We consider the moduli space  $M=M_{2d}$  of K3 surfaces with a 
 polarization of degree $2d$ in characteristic $p$. Let $(X,D)$ be a
polarized K3 surface with a  polarization of degree $2d$. The existence
of this moduli spaces follows from work of Gieseker. We view these moduli spaces
as algebraic stacks.  If the Chern class $c_1(D)$ is not zero in the de Rham
cohomology of $X$ then the moduli space is formally smooth at $[(X,D)]$. 
\par
We shall assume for simplicity that the degree $2d$ of the polarization is prime
to $p$.  Let furthermore $\pi:{\cal X}\tto M_{2d}$ be the universal family of
polarized K3 surfaces over $k$. We set 
$$
       M^{(h)}:=\{ s \in M: h(X_s)\geq h\}.
$$ 
Then, by Artin [A], $M^{(h)}$ is an algebraic subvariety of 
codimension $\leq h-1$ in $M$ for $h=1,\ldots,10$. We shall show that their 
codimension is $h-1$.

The direct image sheaves $R^2\pi_*W_i(O_{\cal X})$ are coherent
sheaves of rings, but not coherent $O_M$-modules. If there would
exist a suitable Grothendieck group of such objects we could
calculate Chern classes by using Theorem 5.1. 
Since we do not know how to do this
we resort to a different method to calculate cycle classes of loci of
given height. 
\par
Let $X_{0}$ be a K3 surface, and assume that the height of the formal
Brauer group $\Phi_{X_{0}}$ is greater than or equal to $h$, i.e.,
$X_{0}$ corresponds to a point in $M^{(h)}$. Then, the Frobenius morphism
is zero on $H^{2}(X, W_{i}(O_{X/S}))$ for $i = 1, \ldots, h - 1$. 
We let $S$ be a formal neighborhood of $M^{(h)}$ at the point, and
we also denote by $\nabla$ the Gauss-Manin connection of 
$H^{2}_{\rm dR}(X/S)$. 
We consider the Hodge filtration 
$0 \subset F^{2} \subset F^{1} \subset H^{2}_{\rm dR}(X/S)$, and construct,
in the same way as in Section 8, a homomorphism 
$$
\Phi_{h} : H^{2}(W_h(O_X)) \tto H^{2}_{\rm dR}(X).
$$
We take a basis $\omega$ of 
$H^{0}(\Omega_{X/S}^{2})$ and take the dual basis $\zeta$ of 
$H^{2}(O_{X/S})$. We take a lifting 
${\tilde \zeta} \in H^{2}_{\rm dR}(X/S)$ of $\zeta$. Then we
have $\langle {\tilde \zeta}, \omega\rangle = 1$. 
Since $R^{n-1}: H^{2}(W_{n}(O_{X/S})) \rightarrow H^{2}(O_{X/S})$
is surjective, we take an element $\alpha \in H^{2}(W_{h}(O_{X/S}))$
such that $R^{h-1}(\alpha) = \zeta$. 
We set 
$$
g_{h} = \langle \Phi_{h}(\alpha), {\omega}\rangle.
$$
Since $\Phi_{h}(\alpha) - g_{h}{\tilde \zeta}$ is orthogonal
to $\omega$, it follows that $\Phi_{h}(\alpha) - g_{h}{\tilde \zeta}$ is contained
in the $F^{1}$-step of the Hodge filtration. Therefore, using  the natural
isomorphism 
$H^{2}_{\rm dR}/F^{1} \cong H^{2}(O_{X})$, we conclude
that
$$
 \phi_{h}(\zeta) = g_{h}{\zeta} \quad {\rm in}~H^{2}(O_{X}),
$$
where $\phi_h$ was defined in section 5. 
This means that the equation $g_{h} = 0$ gives the scheme
theoretic locus of zero of $\phi_{h}$, and by Theorem 8.2, 
the support of the locus in $M^{(h)}$ coincides with 
$M^{(h+1)}$.
\par
\noindent
\proclaim (\chno.1) Proposition. Under the notation and assumptions
made above, the image ${\rm Im}~\Phi_{h}$ is horizontal
with respect to the Gauss-Manin connection. \par
\noindent
{\sl Proof.} It suffices to prove $\nabla (\Phi_{h}(\alpha)) = 0$.
 The element 
$\alpha$ is represented by a cocycle
 $\alpha_{ijk}=(\alpha^{(0)}_{ijk}, \ldots, \alpha^{(h-1)}_{ijk})$ 
with respect to a suitable affine open covering $\{U_{i}\}$ of $X/S$.
Since the Frobenius morphism is zero on $H^{2}(W_{h-1}(O_{X/S}))$,
there exist a cochain $\gamma_{ij}\in \Gamma(U_i\cap U_j,W_{h-1}(O_X/S))$ 
such that $FR(\alpha) = \partial \{ \gamma_{ij}\} = \{
\gamma_{jk}-\gamma_{ik} + \gamma_{ij}\} (\in C_2(W_{h-1}))$. Hence we have
$$ 
F(\alpha)-\partial(\{ (\gamma_{ij}, 0)\})= \{(0,\ldots,0,g_{ijk})\}. \eqno(15)
$$
Put ${\tilde \gamma}_{ij} = (\gamma_{ij}, 0)$, an element in 
$\Gamma(U_{i}\cap U_{j}, W_{h}(O_{X/S}))$.
Then
$\phi_h(\zeta)= \{ g_{ijk}\}$ and
$$
\Phi_h(\alpha)= ( g_{ijk}, -D_h(\tilde{\gamma}_{ij}),0) \in
C_2(O_{X/S})\oplus C_1(\Omega^1_{X/S})\oplus C_0(\Omega^2_{X/S}).
$$
We write this as
$$
\Phi_h(\alpha)=\{( g_{ijk},b_{ij},0)\}
$$
We have to calculate $\nabla(\Phi_{h}(\alpha))$. We use the explicit description
of the Gauss-Manin connection. Katz and Oda define in [K-O] two operators
$$
L_S: C_q(\Omega^p)\to C_q(\Omega^{p+1}), 
\quad L_S((\beta)(i_0,\ldots,i_q)) = d_S^i(\beta(i_0,\ldots,i_q))
$$
and
$$
\lambda: C_q(\Omega^p) \to C_{q+1}(\Omega^p), \quad
\lambda(\beta)(i_0,\ldots,i_{q+1})=(-1)^p(I^{i_0}-I^{i_1})
(\beta(i_1,\ldots,i_{q+1})).
$$
Here we follow the notation of loc.\ cit. The (substitution) operator
$I^{i_0}$ is given by $\sum_{t=1}^p
\hbox{\rm subs}(dx_t \mapsto d_S^{i_0})$ and is zero for $p=0$. In our case
this gives $L_S(g_{ijk})=d_S^i(g_{ijk}) \in C_2(\Omega^1)$,
$\lambda(g_{ijk})=0$ and $L_S(b_{ij})=d_S^i(b_{ij}) \in C_1(\Omega^2)$,
$\lambda(b)(ijk)=-(I^i-I^j)(b_{jk})\in C_2(\Omega^1)$. So we find
$$
\nabla(\Phi_{h}(\alpha))= d_S^ib_{ij} + d_S^ig_{ijk} - I^ib_{jk} +
I^jb_{jk}.\eqno(16)
$$
Here the first term lies in $C_1(\Omega^2)$.

Using $d_S^i$ instead of $d$ we can make an operator $D_{h,S}^i$ similar to
the operator $D_h$ defined by Serre. It is zero on the image of Frobenius
and so the relation (15) gives
$$
-D_{h,S}^i(\partial\{ \gamma_{ij}\})= d_S^i(g_{ijk}).
$$
This says $d_S^i(g_{ijk})= I^i(b_{jk}-b_{ik}+b_{ij})$. Collecting the terms
we get
$$
\nabla(\Phi_{h}(\alpha))= d_S^ib_{ij} + I^i(-b_{ik}+b_{ij}) +I^jb_{jk}.
$$
Put $c_{ij}=-D_h(\gamma_{ij})$. Now note that we have
$$
d_S^i(D_h(\gamma_{ij})) = d(D_{h,S}^i\gamma_{ij}).
$$
Therefore the right hand side of (16) is a boundary in the total complex.
We conclude $\nabla \Phi_{h}(\alpha) = 0$ in 
$\Omega^{1}_{S/k} \otimes H^{2}_{DR}(X/S)$.
$\square$
\smallskip
\bigskip
\def\chno{14}\centerline{\bf \chno. The Tangent Spaces to the Stratification }
\smallskip 
\noindent
We denote by $D_{0}$ the polarization class of $X_{0}$ of degree $2d$
and we shall assume that is not a $p$-th power. Let $M^{(h)}$ be the closed locus
of the moduli space $M=M_{2d}$ of polarized K3 surfaces  given by the condition
height $\geq h$ for $h =1,\ldots, 10$ and $h=\infty$. We now determine the
tangent space of $M^{(h)}$ at the point $x_{0} = (X_{0}, D_{0})$.  We denote by
${\rm Im}~H^{1}(X_{0}, Z_{\ell})$  the image of  $H^{1}(X_{0},
Z_{\ell}\Omega^{1}_{X_{0}})$ in $H^{1}(X_{0}, \Omega^1_{X_{0}})$ induced by
the natural inclusion $Z_{\ell}\Omega^1_{X_{0}} \to \Omega^1_{X_{0}}$.
\par
\noindent
\proclaim (\chno.1) Proposition. Suppose that  $(X_0,D_0)$ represents a point
$x_0$ of $M^{(h)}- M^{(\infty)}$. Then for $1\leq h \leq 10$  the tangent space of
$M^{(h)}$  at $x_{0}$ is in a natural way  isomorphic to 
${\rm Im}~H^{1}(X_{0}, Z_{h-1}) \cap c_{1}(D_{0})^{\perp}$. 
\par
\noindent
{\sl Proof.} Note that by (9.2) the map $H^1(X_0,B_{h-1}) \to  H^1(X_0,
\Omega^1_{X_{0}})$ is injective.  Since we have $H^{1}(X_{0}, B_{h-1}) \subset
c_1(D_{0})^{\perp}$, by Corollary (10.2) and Lemma (9.3), it suffices to prove
that 
$\langle H^{1}(X_{0}, B_{h-1}), c_{1}(D_{0}) \rangle$ is the normal space of 
$M^{(h)}$ at $x_{0}$. We will show this by induction. 
Note that we know $\dim H^{1}(X_{0}, B_{\ell}) = \ell$ for 
$\ell = 0, \ldots, h-1$.
\par
Suppose $h=1$. Then, we have $H^{1}(X_{0}, B_{0}) = 0$,
and by the general theory of moduli spaces the tangent
space of $M^{(1)} = M$ at $x_{0}$ is given by 
$c_{1}(D_{0})^{\perp} \subset H^{1}(X_{0}, \Omega^{1}_{X_{0}})$. 
\par
Now, we assume that the statement holds until $h$.
We use the notation above. Then, by (8.2) $M^{(h+1)}$ is defined by
$g_{h} = \langle \Phi_{h}(\alpha), {\omega}\rangle = 0$ in $M^{(h)}$.  Using
Proposition (13.1), we have $$
\eqalign{
  dg_{h}  & =  \langle \nabla \Phi_{h}(\alpha), \omega \rangle
    +  \langle \Phi_{h}(\alpha), \nabla {\omega}\rangle  \cr
    & = \langle \Phi_{h}(\alpha), \nabla {\omega}\rangle
\cr}
$$
We denote by $m$ (resp. $m_{0}$) the maximal ideal which corresponds
to the point $x_{0}$ in the versal formal moduli space around $x_{0}$
(resp. in the formal moduli around $x_{0}$ in $M^{(h)}$). Then, under  the
natural homomorphism 
$$
    H^{1}(X_{0}, \Omega^{1}_{X_{0}}) \cong m/m^{2} \tto m_{0}/m^{2}_{0}
$$
$- \Phi_{h}(\alpha)(0)$ corresponds to the cotangent vector $g_{h}$ 
by the argument of Ogus [O]. The kernel of this homomorphism is isomorphic
to $\langle H^{1}(X_{0}, B_{h - 1}) , c_{1}(D_{0}) \rangle$ by induction.
We have 
$$
\eqalign{
- \Phi_{h}(\alpha)(0)&=  -\{D_{h}(\tilde \gamma_{ij})\}   \cr
  &= -\{ \sum_{m=0}^{h-1} (\gamma_{ij}^{(m)})^{p^{h-m}-1}d\log
\gamma_{ij}^{(m)} \} \cr                    }
$$
and
$D_{h} : H^{2}(W_{h}(O_{X_{0}})/FW_{h}(O_{X_{0}})) \to 
 H^{1}(X_{0}, \Omega^{1}_{X_{0}})$ is injective 
by Corollary (9.2).
Since $\Phi(\alpha)(0)$ lies in $H^1(X_0,B_{h})$ but not in $H^1(X_0,B_{h-1})$,
we conclude that $g_{h} \notin m^{2}_{0}$.
By induction we thus see that the tangent space to $M^{(h+1)}$ can be
identified with $H^1(X_0,Z_{h})\cap c_1(D_0)^{\bot}$. $\square$
\smallskip
This argument does not work for $h=\infty$, but can be made to work
for the supersingular points for which the subspace 
$\langle{\rm Im}(H^1(B_h)), c_1(D)\rangle$ of $H^1(\Omega^1)$  has
dimension $h$. In Section 12 we gave conditions for this. Under the assumption
that $\rho=B_2$ this is the case if the Artin invariant $\sigma_0$ of a
supersingular K3 surface satisfies $\sigma_0 > h$. We thus find:
\par
\proclaim (\chno.2) Theorem. For $h=1,\ldots, 10$ the open stratum $M^{(h)}$ if
not empty  is purely of dimension $(20-h)$ and nonsingular at any
point of the stratum $M^{(h)}$ where  the subspace $\langle{\rm Im}(H^1(B_{h-1})),
c_1(D_0)\rangle$ of $H^1(\Omega^1_{X_0})$  has dimension $h$.  In particular, it
is non-singular at non-supersingular points and assuming the Artin conjecture at
all supersingular points with Artin invariant $\sigma_0 \geq h$ and $c_1(D_0)
\not\in {\rm Im}(H^1(B_h))$.
\par
\smallskip
We refer here to a forthcoming preprint of Ogus for a description of the
singularities of the strata. Ogus proved in [O, Prop.\ 2.6] that for $p\neq 2$
the stratum $M^{(2)}$ has a quadratic singularity at the superspecial points. A
variation of his argumentthere shows that at a point with Artin
invariant $\sigma_0=h-1$ the singular locus has multiplicity $2$. In particular
the stratum $M^{(11)}$ has multiplicity $2$ at points with $\sigma_0=10$, cf.\ his
forthcoming preprint and the discussion in the next section.
\par
\bigskip
\def\chno{15}
\noindent
\centerline{\bf \chno. The Loci of K3 Surfaces of Given Height}
\smallskip
\noindent
We now come to the description of the cycle classes of the
strata defined by the height.  Let $M^{(h)}$ be the closed stratum of the moduli
space $M=M_{2d}$ where the height of the formal group $\Phi_X$ is at least $h$
with the convention that
$M^{(11)}=M^{(\infty)}$. For simplicity we shall assume that $p$ does not divide
$2d$. By our characterization of $h$ these strata can be given a natural scheme
structure and these are reduced for $h\neq
\infty$ by our results in Section 14. It is known by Artin  that the strata
$M^{(h)}$ for $h=1,\ldots,11$ have codimension $\leq h -1$ in $M_{2d}$, see [A]. 
\par
Define a line  bundle $V$ on $M$ by  $V=\pi_*(\Omega^2_{{\cal X}/M})$ and let the first Chern class be $v$. 
\noindent
\proclaim (\chno.1) Theorem. Let $M=M_{2d}$ be the moduli stack of polarized K3
surfaces over $k$ with a polarization of degree $2d$ prime to $p$.  Then
for $h=1,\ldots, 10, 11$ the scheme-theoretic locus $M^{(h)}$ of surfaces
with height $\geq h$, if not empty, is of codimension $h-1$ and for $h\neq 11$ it 
is a local complete intersection. The class of $M^{(h)}$ in the Chow group
$CH_{\Q}^{h-1}(M)$ is given by
$$
(p-1)(p^2-1)\ldots(p^{h-1}-1)v^{h-1}.
$$
\par
\noindent
{\sl Proof.} We prove this by induction. Let $M$ be the moduli space
of polarized K3 surfaces of degree $2d$ as above. We know that the generic K3
surface has height $1$, and so  for $h=1$ the formula is correct. The
codimension of $M^{(h)}$ is $\leq h-1$ for $1\leq h\leq 10$ as follows from
(5.7). For $h=2$ the locus $M^{(2)}$ is the non-ordinary locus. This locus is
characterized by the fact that the Frobenius map $H^2(X,O_X) \to
H^2(X,O_X)$ vanishes. This is a $p$-linear map and the corresponding $O_M$-linear
map is $(R^2\pi_*O_{\cal X})^{(p)} \to R^2\pi_*O_{\cal X}$ with associate
cycle class $(p-1)v$. Locally, at a point of $M^{(2)}$ an equation is
given by $g_1=0$, see the proof of (14.1) and $dg_1\neq 0$. So if $M^{(2)}$ is not
empty then it is purely $18$-dimensional.
\par
Suppose now that the class of $M^{(h-1)}$ is given by the
class in the formula. By Proposition (5.7) the locus in $M^{(h-1)}$
where the height increases is given by the vanishing of the map
$\phi_{h-1}: (R^2\pi_*O_{\cal X})^{(p^{h-1})} \to R^2\pi_*(O_{\cal
X})$, equivalently, by the vanishing of a section of $V^{p^{h-1}-1}$. By (14.1)
it follows that for a local equation $g_h=0$ we have $dg_h\neq 0$. Hence the 
locus is reduced for $h\neq\infty$ and  the class
on $M^{(h-1)}$ is given by $(p^{h-1}-1)v$.  
\par
Let $j_h: M^{(h)} \to M^{(h-1)}$ and $j: M^{(h-1)} \to M$ be the natural
inclusions. Then the class of the locus $M^{(h)}$ in $CH_{\Q}^{h-1}(M)$
is given by 
$$
{j}_* {{j}_h}_* [M^{(h)} ]= {j}_* ([M^{(h-1)}]\cdot
j_{h-1}^* (p^{h-1}-1)v)=(p^{h-1}-1) v \cdot
({j}_*[M^{(h-1)}])
$$
by the projection formula. $\square$
\par
\smallskip
The locus $M^{(11)}$ comes with a multiplicity in the formula because of (14.2).
For $p\neq 2$ the multiplicity is $2$.  It makes sense to call
the reduced locus
$M^{(11)}_{\rm red}$ the supersingular locus. 
\smallskip
\noindent
{\bf (\chno.2) Remark.} In [G] a formula for the class of the supersingular locus
on the moduli space of principally polarized abelian surfaces was given. 
Comparison with Kummer surfaces shows that this is compatible  with 
multiplicity $2$  along  the supersingular locus,  cf.\  [G-K].
\bigskip
We shall now assume that the line bundle $V=\pi_*(\Omega^2_{{\cal X}/M})$ is 
ample on the moduli space. It is known by the theory of Baily and Borel (see
[B-B]) that $V$ is ample on the moduli spaces in characteristic $0$; indeed,
modular forms of sufficiently high weight define an embedding.
\smallskip
\noindent
\proclaim (\chno.3) Theorem. Suppose that the class $v$ is ample. Let
${\cal X} \to S$ with $S$ complete be a proper smooth family of  polarized K3
surfaces with constant $h \neq \infty$. Then this family is isotrivial.
\par
\noindent
{\sl Proof.} It follows from the preceding theorem that the strata
$S^{(h)} -S^{(h +1)}$ where the height is constant are quasi-affine
for $h=1,\ldots,10$.
\par 
\smallskip
We do not know whether the class $v$ is ample on the moduli spaces
$M_{2d}$, but we expect it to be so.    
\smallskip
Suppose that there exists a good Baily-Borel compactification. By this
we mean that there exists a projective variety (stack) ${\overline
M}_{2d}$ containing $M_{2d}$ such that ${\overline M}_{2d}-M_{2d}$ is
$1$-dimensional and consists of a configuration of elliptic modular
curves. This is the case in characteristic zero, cf.\ Kondo [Ko]. Then it
follows from our theorem that a family of K3 surfaces with $h\geq 3$
does not degenerate. Indeed, it follows from our formula that a class
of the form  $v^m$ with $m\geq 3$ has zero intersection with the
`boundary components'. This implies that for each boundary component
the locus with $h\geq 3$ either has empty intersection with this
boundary component or contains it. The boundary components form a
connected set and the generic point of each component crresponds to  a degenerate
K3 surface corresponding to an ordinary elliptic curve. For the degenerate
surfaces the  height is $1$ or $2$. Compare the discussion in [R-Z-Sh].
\bigskip
\bigskip
\noindent
\def\chno{16}\centerline{\bf \chno. An Extension for Other Varieties}
\smallskip
\noindent
Though the theorem in Section 5 was formulated for K3 surfaces
it holds for a more general class of surfaces. 
\proclaim (\chno.1) Theorem. Suppose that $X$ is a smooth algebraic surface
such that \item{i)} ${\rm Pic}^0(X)$ is reduced, 
 \item{ii)} $\dim H^2(X,O_X)=1$. 
\smallskip 
Then  $\Phi^2$ is represented by a formal group of dimension $1$ and its
height satisfies $h(\Phi_X)
\geq i+1$  if and only if the Frobenius map $ F$ on  $H^2(X,W_i(O_X))$ is
the zero map.
\par
\proclaim (\chno.2) Corollary. For such a surface we have the following 
characterization of the height:
$$
h(\Phi_X) = \min \{ i \geq 1 \colon [F: H^2(W_i(O_X)) \to
H^2(W_i(O_X))] \neq 0 \}.
$$\par
\smallskip
\noindent
{\sl Proof.} The proof is analogous to the proof given for K3 surfaces.  Instead
of the vanishing of $H^1(X,O_X)$ one uses the vanishing of the Bockstein
operators. Recall that $H^1(X,W_n(O_X))$ is the subgroup of
$k[\epsilon]/(\epsilon^{n+1})$-valued points of the connected component of the
Picard scheme $P$ at the origin, cf.\ [Mu]. A $k[\epsilon ]/\epsilon^2$-valued
point (tangent vector) is tangent to $P_{\rm red}$ at the origin if and only if
it can be lifted to $k[\epsilon]/(\epsilon^n)$-valued point for all $n$. That is,
these correspond precisely to the elements of $H^1(X,O_X)$ that can be lifted to
$H^1(X,W_n(O_X))$ for all $n$. So if $P=P_{\rm red}$ then all elements of
$H^1(X,O_X)$ can be lifted and this implies the analogues of Lemmas (4.2) and
(4.5) that we need.
\smallskip
\noindent
{\bf \chno.3) Example.} 1) An abelian surface satisfies the assumptions. 
2) A surface of general type with $H^1(O_X)=0$ and $p_g=1$. Examples of such
surfaces are surfaces with $K^2=p_g=1$. These have $h^1(X,O_X)=0$ and are
resolutions of surfaces of type (6,6) in weighted projective space
$\P(1,2,2,3,3)$, cf.\ [C].
\smallskip
A nonsingular complete algebraic variety $X$ of dimension $n$ is called
a Calabi-Yau variety if the canonical invertible sheaf $\omega_X$
is trivial and $H^{i}(X, O_{X}) = 0$ for $1 \leq i \leq n - 1$. By
a criterion of Artin-Mazur [A-M], the Artin-Mazur formal group $\Phi^{n}$
is pro-representable by a one-dimensional formal Lie group for such a variety.
In the same way as in Section 5, we have also a characterization
of the height of the formal group~$\Phi^{n}$.
\par
\proclaim (\chno.4) Proposition. For a Calabi-Yau variety $X$ of dimension $n$ we
have the following characterization of the height:
$$
h(\Phi^{n}_{X}) = \min \{ i \geq 1 \colon [F: H^{n}(W_i(O_{X})) \to
H^{n}(W_{i}(O_X))] \neq 0 \}.
$$
\par
\smallskip
\noindent\bigskip
\centerline{\bf References}
\smallskip
\noindent
[A] M.\ Artin: Supersingular K3 surfaces.  {\sl Ann.\ Scient.\ Ec.\ Norm.\
Sup.\ \bf 7} (1974), p.\ 543--568.
\smallskip
\noindent
[A-M] M.\ Artin, B. Mazur: Formal groups arising from algebraic
varieties. {\sl Ann.\ Scient.\ Ec.\ Norm.\ Sup.\
\bf 10} (1977), p.\ 87--132.
\smallskip
\noindent
[B-B] W.L.\ Baily, Jr.; A.\ Borel: Compactification of arithmetic
quotients of bounded symmetric domains. {\sl Ann.\ of Math.\ (2) \bf
84} (1966) p.\ 442--528. 
\smallskip
\noindent
[C] F.\ Catanese:   Surfaces with $K^{2}=p_{g}=1$ and their period mapping.
Algebraic geometry (Proc. Summer Meeting, Univ. Copenhagen, Copenhagen, 1978),
 {\sl Lecture Notes in Math. \bf 732,} p.\ 1-29, Springer Verlag 1979.
\smallskip
\noindent
[D] P.\ Deligne: Rel\`evement des surfaces K3 en
caract\'eristique nulle. Prepared for publication by Luc Illusie.
{\sl Lecture Notes in Math.}, 868, Algebraic surfaces (Orsay,
1976--78), p.\ 58--79, Springer, Berlin-New York, 1981. 
\smallskip
\noindent
[G] G.\ van der Geer: Cycles on the moduli space of abelian varieties.
In: {\sl Moduli of Curves and Abelian Varieties.} The Dutch
Intercity Seminar on Moduli. Eds.  C.\ Faber and E.\ Looijenga. 
Aspects of Mathematics, Vieweg 1999.
\smallskip
\noindent
[G-K] G.\ van der Geer, T.\ Katsura: Formal Brauer groups and the moduli
of abelian surfaces. {\tt math.AG/9912169} 
\smallskip
\noindent
[H] Hazewinkel: Formal Groups and Applications. Boston Orlando, 
Academic Press, 1978.
\smallskip
\noindent
[Il] L.\ Illusie: Complexe de de Rham-Witt et cohomologie cristalline. 
{\sl Ann.\ Sci.\ ENS \bf 12} (1979),p.\ 501--661.
\smallskip
\noindent
[K-O] N.\ Katz, T.\ Oda: On the differential of de Rham cohomology 
classes with respect to parameters. {\sl J. Kyoto Univ. \bf 8} 
(1968), p.\ 199-213.
\smallskip
\noindent
[Ko] S.\ Kondo: On the Kodaira dimension of the moduli spaces of K3 
surfaces.
{\sl Compositio Math., \bf 89}, (1993), p.\ 251--299.
\smallskip
\noindent
[Ma] Yu.\ I.\ Manin : The theory of commutative formal groups over
fields of finite characteristic. {\sl Russian Math.\ Surv.\ \bf 18}, (1963), p.\
1-80 (= {\sl Usp.\ Mat.\ Nauk.\ \bf 18}, (1963), p.\ 3--90.
\smallskip
\noindent
[Mu] D.\ Mumford: Lectures on Curves on an Algebraic Surface.  Annals of
Math.\ Studies 59. Princeton University Press 1966.
\smallskip
\noindent
[O] A.\ Ogus: Supersingular K3 crystals. {\sl Ast\'erisque \bf 64} (1979), 
p.\ 3-86.
\smallskip
\noindent
[R-Sh] A.N.\ Rudakov, I.\ R.\ Shafarevich:  Surfaces of type K3 over fields of
finite characteristic. {J. Soviet Math.} (1983), p.\ 1476-1533.
\smallskip
\noindent
[R-Z-Sh] A.N.\ Rudakov, T.\ Zink, I.R.\  Shafarevich.   The effect of height
on degenerations of algebraic $K3$ surfaces. (Russian) {\sl Izv.\ Akad.\ Nauk
SSSR Ser.\ Mat.\ 46} (1982), p.\ 117--134, 192.
\smallskip
\noindent
[Sc] F.\ Scattone: On the compactification of moduli spaces for  algebraic K3
surfaces. {\sl Memoirs of the A.M.S., \bf 70}, No347 (1987).
\smallskip
\noindent
[S] J-P. Serre: Sur la topologie des vari\'et\'es alg\'ebriques en
caract\'eristique $p$. {\sl Symposion Internacional de topologia
algebraica} 1958, p.\ 24--53.
\smallskip
\noindent
[Sh] T.\ Shioda: Supersingular K3 surfaces.  
Algebraic Geometry  (Proc. Summer Meeting, Univ. Copenhagen, 
Copenhagen, 1978),
{\sl Lecture Notes in Math., {\bf 732}}, p.\ 564-591, Springer-Verlag 1979. 
\bigskip
\settabs4 \columns
\+G.\ van der Geer  &&T.\ Katsura\cr
\+Faculteit
Wiskunde en Informatica &&Graduate School of Mathematical Sciences\cr
\+Universiteit van
Amsterdam && University of Tokyo\cr
\+Plantage Muidergracht 24&&3-8-1 Komaba, Meguro-ku \cr
\+1018 TV Amsterdam
&&Tokyo 153-8914 \cr
\+The Netherlands &&Japan \cr
\+{\tt geer@wins.uva.nl} &&{\tt tkatsura@ms.u-tokyo.ac.jp} \cr

\end